\documentclass[11pt]{article}
\usepackage{amsfonts}
\usepackage{amsfonts}
\usepackage{amssymb,amsmath}
\usepackage{array}
\usepackage{mathrsfs}
\usepackage{mathrsfs}
\usepackage{amsfonts}
\usepackage{fancyhdr}\usepackage{framed}
\usepackage{latexsym}
\usepackage{bbding}
\usepackage{wasysym}
\usepackage{cite}
\usepackage{multicol,graphics}
\usepackage{amssymb}
\usepackage{graphicx}
\usepackage{color}
\usepackage{enumitem}
\usepackage{indentfirst} % 确保首行缩进
 \allowdisplaybreaks

\catcode`!=11
\let\!int\int \def\int{\displaystyle\!int}
\let\!lim\lim \def\lim{\displaystyle\!lim}
\let\!sum\sum \def\sum{\displaystyle\!sum}
\let\!sup\sup \def\sup{\displaystyle\!sup}
\let\!inf\inf \def\inf{\displaystyle\!inf}
\let\!cap\cap \def\cap{\displaystyle\!cap}
\let\!max\max \def\max{\displaystyle\!max}
\let\!min\min \def\min{\displaystyle\!min}

\catcode`!=12

\let\oldsection\section
\renewcommand\section{\setcounter{equation}{0}\oldsection}

\allowdisplaybreaks
\def\pf{\it{Proof.}\rm\quad}

\newcommand{\mdx}{\;\mathrm{d}x}
\newcommand{\md}{\;\mathrm{d}}

\newcommand\divg{{\rm{div}}}

\newtheorem{thm}{Theorem}[section]
\newtheorem{pro}{Proposition}[section]

\newtheorem{lem}{Lemma}[section]
\newtheorem{re}{Remark}[section]
\begin{document}
\title{Global regularity and incompressible limit of 2D compressible Navier-Stokes equations with large bulk viscosity\thanks{This work was partially supported by the National Natural Science Foundation of China (Nos.  12471226, 12131007, 12071390, 12226344, 12226347), and the Natural Science Foundation of Liaoning Province (Grant No. 2023-MS-137)}}
\author{Shengquan Liu$^1$, \quad Jianwen Zhang$^{2,3}$\thanks{Emails: shquanliu@163.com (S. Q. Liu),   jwzhang@xmu.edu.cn (J. W. Zhang, Corresponding author) }\\
\small $^1$ School of Mathematics and Statistics, Liaoning University, Shenyang 110036, P. R. China\\ 
\small $^2$ School of Science, Jimei University, Xiamen 361021, P. R. China\\
\small $^3$ School of Mathematical Sciences, Xiamen University, Xiamen 361005, P. R. China}
\date{}
\maketitle \noindent{\bf Abstract.} In this paper, we study the global regularity of large solutions with vacuum to the two-dimensional compressible Navier-Stokes equations  on $\mathbb{T}^{2}=\mathbb{R}^{2}/\mathbb{Z}^{2}$, when the volume (bulk) viscosity coefficient $\nu$ is sufficiently large. It firstly fixes a flaw in \cite[Proposition 3.3]{Danchin2023}, which concerns the $\nu$-independent global $t$-weighted estimates of the solutions. Amending the proof requires non-trivially mathematical analysis. As a by-product, the incompressible limit with an explicit rate of convergence is shown, when the volume viscosity tends to infinity. In contrast to \cite[Theorem 1.3]{Danchin2019} and \cite[Corollary 1.1]{DM2017}  where vacuum was excluded, the convergence rate of the incompressible limit is obtained for the global solutions with vacuum, based on some $t$-growth and singular $t$-weighted estimates.
\\[2mm]
\\[2mm]
\section{Introduction}\label{sec1}

This paper is concerned with the baratropic  compressible Naveir-Stokes equations on a unit torus $\mathbb{T}^{2}=\mathbb{R}^{2}/\mathbb{Z}^{2}$:
\begin{align}
\begin{cases}\label{1.1}
 \partial_{t} \rho+\divg(\rho u)=0,\\[1mm]
  \partial_{t}(\rho u)+\divg(\rho u \otimes u)+\nabla P(\rho )=\mu\Delta u+
(\mu+\lambda)\nabla\divg u,
\end{cases}
\end{align}
which is supplemented with the initial data
\begin{equation}\label{1.2}
(\rho, u)|_{t=0}=(\rho_{0}, u_{0})(x),\quad x\in\mathbb{T}^2.
\end{equation}
Here, the unknown functions $\rho=\rho(x,t)\geq0$ and $u=(u_1,u_2)(x,t)$ are the fluid 
density and the velocity, respectively. The pressure $P=P(\rho)$ is a nonnegative $C^1$-function of the density and satisfies $P'(\rho)\geq0$ for all $\rho\geq0$. For simplicity, we focus on the physically relevant case that the pressure obeys  the $\gamma$-law equation of state,
\begin{equation}\label{1.3}
P(\rho)=A\rho^\gamma\quad {\text{with}}\quad \gamma\geq 1\quad{\text {and}}\quad A>0.
\end{equation}
The real numbers $\mu$ and $\lambda$ are the shear and bulk viscosity coefficients, respectively, and  satisfy 
\begin{equation}\label{1.4}
\mu>0\quad{\text {and}}\quad  \nu:=2\mu+\lambda> 0.
\end{equation}

The mathematical theory of compressible Navier-Stokes equations has been  extensively studied by many people, due to its physical importance and mathematical challenges, see, for example, \cite{Des1997,Feireisl,Ho1995,HLX2012,Lions,MN1980,Nash,Solo1980} and among others. However, the question of regularity and uniqueness of global solutions with generally large data remains open. In \cite{Kazhikhov1995,HL2016,Fan2022}, the authors studied the global well-posedness of large solutions of 2D compressible Navier-Stokes equations with density-dependent bulk viscosity. Recently,   Danchin-Mucha \cite{DM2017,Danchin2019,Danchin2023} obtained the global solutions of (\ref{1.1}) for a class of large data. In particular,  Danchin-Mucha \cite{Danchin2023}  proved  the existence of global solutions to the problem (\ref{1.1})--(\ref{1.4}) with ripped density (see \cite[Theorem 2.1]{Danchin2023}),  provided  the  volume (bulk) viscosity coefficient $\nu$ is large enough. It was assumed  in \cite{Danchin2023} that the initial velocity $u_0\in H^1$ and  the initial density $0\leq\rho_0\in L^\infty$, and hence, the solutions may contain  vacuum states. Moreover, the uniqueness of global solutions was shown to hold (cf. \cite[Theorem 2.2]{Danchin2023}) in the case that the pressure $P(\rho)$ depends linearly on the density (i.e. the isothermal system with $\gamma=1$ and $P(\rho)=A\rho$). As a by-product of the uniform-in-$\nu$ estimates, the incompressible limit (without convergence rates) from the compressible Navier-Stokes equations (\ref{1.1}) to the inhomogeneous incompressible Navier-Stokes equations was also justified in \cite[Theorem 2.3]{Danchin2023}, as the bulk viscosity tends to infinity. We also refer to \cite{DM2017,Danchin2019} for analogous results with strictly positive density.

Due to the lack of the estimate $\nabla u\in L^1(0,T;L^\infty)$ in \cite{Danchin2023}, it seems difficult to reformulate system (\ref{1.1})  in Lagrangian coordinates and to prove the uniqueness in the same way as that in \cite{DFP2020}. However, Danchin-Mucha \cite{Danchin2023} succeeded in proving the uniqueness of the solutions for the isothermal Navier-Stokes equations (\ref{1.1}) with $P(\rho)=A\rho$ in an analogous manner as that in \cite{Ho2006} by  bounding the difference of the density in $H^{-1}$-norm and making use of a backward parabolic system to deal with the difference of velocity, provided the solutions satisfy  $(\divg u,\nabla \mathcal{P}u)\in L^r(0,T; L^\infty)$ (and hence, $\nabla u\in L^r(0,T; {\text{BMO}})$) for some $1<r<2$.  In  fact, the estimates of  $(\divg u,\nabla\mathcal{P}u)\in L^r(0,T; L^\infty)$ come from the following $t$-weighted estimates (cf. \cite[Proposition 3.3]{Danchin2023}) and the regularity theory of elliptic system (see \cite[Corollary 3.1]{Danchin2023}). Note that the $t$-weighted estimates are essentially needed to circumvent the difficulty induced by the lower regularity of initial data. 

\begin{pro}[{\cite[Propsition 3.3]{Danchin2023}}] \label{pro1.1}Assume that $(\rho,u)$ is a smooth solution of the problem \eqref{1.1}--\eqref{1.4} on $\mathbb{T}^2\times[0,T]$ with $0<T<\infty$. Then, there exists a positive constant $C(T)$, independent of $\nu=2\mu+\lambda$, such that for sufficiently large $\nu\geq1$,
\begin{equation}\label{1.5}
\sup_{0\leq t\leq T}t\int_{\mathbb{T}^2} \rho|\dot u |^2\mdx+\int_0^Tt\int_{\mathbb{T}^2}\left(\mu|\nabla \mathcal{P}\dot u|^2+\nu|\divg \dot u|^2\right)\mdx\leq C(T),
\end{equation}
where $\dot f:=f_t+u\cdot\nabla f$ is the material derivative of $f$ and $\mathcal{P}:=\mathbb{I}+\nabla(-\Delta)^{-1}\divg$ denotes the $L^2$-projector onto the set of solenoidal vector fields.
\end{pro}

Let us make some comments on the proof of (\ref{1.5}) in \cite{Danchin2023}. Indeed, the largeness of $\nu$ was determined in  \cite[Propositions 3.1 and 3.2]{Danchin2023} and was used to derive the global $H^1$-estimate of the velocity and an upper bound of the density, both of which are uniform in $\nu$ and independent of $t$. Analogously to that in \cite{Ho1995},  to prove (\ref{1.5}), Danchin and Mucha \cite{Danchin2023} operated the material derivative to both sides of the momentum equation (\ref{1.1})$_2$ and tested it by $t\dot u$. Unfortunately, {\it an error} occurs in the handling of the following item:
\begin{equation}\label{1.6}
-\nu\int_{\mathbb{T}^2} \frac{D}{Dt}\nabla \divg u \cdot (t\dot u)\mdx\quad{\text{with}}\quad \frac{D}{Dt}f=\dot f=f_t+u\cdot\nabla f.
\end{equation}
To deal with the term in (\ref{1.6}), Danchin-Mucha  claimed in Step 3 of   \cite[Proposition 3.3]{Danchin2023} that
\begin{equation}\label{1.7}
\begin{aligned}
- \frac{D}{Dt}\nabla \divg u &= -\nabla \frac{D}{Dt}\divg u+\nabla u\cdot\nabla \divg u\\
&=-\nabla\divg\dot u+\nabla(\text{tr}(\nabla u\cdot\nabla \mathcal{Q}u))+\nabla u\cdot\nabla \divg u,
\end{aligned}
\end{equation}
where $\mathcal{Q}u:=-\nabla(-\Delta)^{-1}\divg u$ denotes the  potential part of the velocity. Note that $u=\mathcal{P}u+\mathcal{Q}u$. However, a direct calculation gives
\begin{equation}\label{1.8}
\begin{aligned}
- \frac{D}{Dt}\nabla \divg u &= -\nabla \frac{D}{Dt}\divg u+\nabla u\cdot\nabla \divg u\\
&= -\nabla\divg\dot u+\nabla\left(\sum_{i,j=1}^2\partial_iu^j\partial_ju^i\right)+\nabla u\cdot\nabla \divg u\\
&=-\nabla\divg\dot u+\nabla(\text{tr}(\nabla u\cdot\nabla u))+\nabla u\cdot\nabla \divg u.
\end{aligned}
\end{equation}
Since these two terms $\nabla(\text{tr}(\nabla u\cdot\nabla \mathcal{Q}u))$ and $\nabla(\text{tr}(\nabla u\cdot\nabla u))$ are entirely different, the derivation of (\ref{1.7}) is incorrect. Indeed, it is obvious that  $\nabla(\text{tr}(\nabla u\cdot\nabla \mathcal{P}u))\neq 0$, though  $\divg \mathcal{P} u=0$ and $\divg u=\divg \mathcal{Q}u$. As it will be seen later, the treatment of the term $\nabla(\text{tr}(\nabla u\cdot\nabla u))$ (see Step 3 of Lemma \ref{lem3.2}) is significantly more sophisticated than that of $\nabla(\text{tr}(\nabla u\cdot\nabla \mathcal{Q}u))$, when one attempts to establish a $\nu$-independent global $t$-weighted estimate analogous to  the one in (\ref{1.5}). 

To clarify the inherent difficulty, we present a detailed explanation of the proof in \cite{Danchin2023}. Let $\widetilde G=\nu\divg u-(P-\overline P)$ be the so-called ``effective viscous flux", where $\overline P$ is the mean of $P(\rho)$ on $\mathbb{T}^2$. Then,  based on integration by parts and   \cite[Propositions 3.1 and 3.2]{Danchin2023}, the  term $ \nu\langle t \dot u, \nabla(\text{tr}(\nabla u\cdot\nabla \mathcal{Q}u))\rangle$ can be easily bounded by (see \cite[(3.64)]{Danchin2023})
\begin{align*}
&\nu\left|\int_0^T t\langle \dot u, \nabla(\text{tr}(\nabla u\cdot\nabla \mathcal{Q}u))\rangle\md t\right| \leq  C \nu\int_0^T\int_{\mathbb{T}^2} t|\nabla u||\nabla \mathcal{Q}u||\divg \dot u|\mdx\md t\\
&\quad \leq C\|\nabla u\|_{L^4(0,T;L^4)}\left(\|G\|_{L^4(0,T;L^4)}+\|P-\overline P\|_{L^4(0,T;L^4)}\right)\|\sqrt t\divg \dot u\|_{L^2(0,T;L^2)}\\
&\quad \leq C\nu^{\frac{1}{4}}\|\sqrt t\divg \dot u\|_{L^2(0,T;L^2)}.
\end{align*}
The key role in the above analysis is to  replace $\nu\divg u$ by $\widetilde G+(P-\overline P)$, since the global estimates in \cite[Propositions 3.1 and 3.2]{Danchin2023} ensure that the quantities  $\|\nabla u\|_{L^4(0,T;L^4)}$, $\nu^{-\frac{1}{4}}\|P-\overline P\|_{L^4(0,T;L^4)}$ and $\nu^{-\frac{1}{4}}\|G\|_{L^4(0,T;L^4)}$ are uniformly bounded in $\nu$. Clearly, this idea is intrinsically incapable of handling the {\it genuine term} $\nu\langle t\dot u, \nabla(\text{tr}(\nabla u\cdot\nabla u))\rangle$, due to the fact that the term $\nabla(\text{tr}(\nabla u\cdot\nabla\mathcal{P} u))$ is non-negligible and must be taken into account. Unfortunately, when $\nu$ is large, the quantity $\nu^{\frac{1}{4}}\|\nabla \mathcal{P}u\|_{L^4(0,T;L^4)}$ (and thus,  $\nu^{\frac{1}{4}}\|\nabla u\|_{L^4(0,T;L^4)}$) cannot be uniformly bounded in $\nu$. Note that $\mathcal{P}u$ is the divergence-free part of the velocity. The lack of these uniform estimates precludes a fully rigorous  proof of the $\nu$-independent estimates stated in (\ref{1.5}). 

Based upon the above observation, the first purpose of this paper is to develop some new ideas to prove an analogous estimate as that in (\ref{1.5}), which also enables us to derive the bounds of  $(\nabla\mathcal{P}u, \divg u)\in L^r(0,T;L^\infty)$ for some $1<r<\infty$. Indeed, instead of bounding  $\nu^{\frac{1}{2}}\|\divg \dot u\|_{L^2(0,T;L^2)}$, it will be shown that $\nu^{\frac{1}{2}}\|D_t(\divg u)\|_{L^2(0,T;L^2)}$ with $D_t=\frac{D}{Dt}$ being the material derivative is uniformly bounded in $\nu$ for any $0<T<\infty$. This, together with the regularity theory of elliptic system and the interpolation arguments, also yields the desired bounds of  $(\nabla\mathcal{P}u, \divg u)\in L^r(0,T;L^\infty)$ for some $1<r<\infty$. As a result, the uniqueness of the solutions to the compressible, viscous and isothermal  Navier-Stokes equations (\ref{1.1}) with $P(\rho)=A\rho$ can be shown to hold. 

To state our main result, we first impose the initial conditions on $(\rho_0,u_0)$. Analogously to that in \cite{Danchin2023}, we assume that there exists a positive number $K>0$ such that
\begin{equation}\label{1.9}
0\leq\rho_0\in L^\infty(\mathbb{T}^2),\quad \int_{\mathbb{T}^2} \rho_0(x)\mdx=1,\quad u_0\in H^1({\mathbb{T}^2}),\quad \|\divg u_0\|_{L^2}\leq K\nu^{-\frac{1}{2}}.
\end{equation} 
Let $F$ and $\dot{u}$ be  the ``effective viscous flux"  and the material derivative defined by
\begin{equation}\label{1.10}
F:=\nu\divg u-(P-\overline P)\quad {\text{with}}\quad \overline P:=\int_{\mathbb{T}^2} P(\rho)\mdx,\quad \dot u:=u_t+u\cdot\nabla u.
\end{equation}

Our first result concerning the  global $\nu$-independent and $t$-weighted estimates of $(\rho,u)$ can now be formulated in the following proposition.

\begin{pro}\label{pro1.2} Let the initial conditions of \eqref{1.9}  be in force. Assume that
$(\rho,u)$ is a smooth solution of the problem \eqref{1.1}--\eqref{1.4} on $\mathbb{T}^2\times[0,T]$ with $0<T<\infty$. Then, there exists a positive constant $C(T)$, independent of $\nu=2\mu+\lambda$, such that for sufficiently large $\nu\geq1$,
\begin{equation}\label{1.11}
\sup_{0\leq t\leq T}t\int_{\mathbb{T}^2} \rho|\dot u |^2\mdx+\int_0^Tt\int_{\mathbb{T}^2}\left(\mu|\nabla^\bot\cdot \dot u|^2+\nu|\dot V|^2\right)\mdx\leq C(T),
\end{equation}
where $V:=\divg u$, $\dot V=V_t+u\cdot\nabla V$, and  $\nabla^\bot\cdot u=\partial_1u_2-\partial_2u_1$ with $\nabla^\bot:=(-\partial_2,\partial_1)$ being the 2D vorticity operator.  Moreover, there also exists a positive constant $C(T)$, independent of $\nu=2\mu+\lambda$, such that  for any $2\leq q <\infty$ and $0<\varepsilon<1$,
\begin{equation}\label{1.12}
\begin{aligned}
&\int_0^T \left(\|\rho\dot{u}\|_{L^q}^{2-\varepsilon}+\|\nabla F\|_{L^q}^{2-\varepsilon}+\|\nabla(\nabla^\bot\cdot  u)\|_{L^q}^{2-\varepsilon}\right)\md t\\
&\quad+\int_0^T\left(\|F\|_{L^\infty}^{2-\varepsilon}+ \|\nabla^\bot\cdot u\|_{L^\infty}^{2-\varepsilon}+\nu^{2-\varepsilon}\|\divg u\|_{L^\infty}^{2-\varepsilon}\right)\md t\leq C(T).
\end{aligned}
\end{equation}
\end{pro}

With the help of Proposition \ref{1.2} and \cite[Propositions 3.1, 3.2]{Danchin2023}, we can prove the global existence theorem of the problem (\ref{1.1})--(\ref{1.4}) in the exactly same way as that in \cite{Danchin2023}.

\begin{thm}\label{thm1.1}
Assume that the initial data $(\rho_0,u_0)$ satisfies  \eqref{1.9}. There exists a positive number $\nu_0$, depending only on $\mu$, $A$, $\gamma$, $K$ and the initial norms of $(\rho_0,u_0)$, such that if $\nu\geq \nu_0$, then the problem \eqref{1.1}--\eqref{1.4} admits a global-in-time solution $(\rho,u)$ on $\mathbb{T}^2\times(0,\infty)$, satisfying 
\begin{equation}\label{1.13}
\begin{cases}
0\leq \rho\in L^\infty(0,\infty; L^\infty(\mathbb{T}^2))\cap C([0,\infty); L^p(\mathbb{T}^2)),\quad \sqrt\rho u\in C([0,\infty); L^2(\mathbb{T}^2)),\\[1mm]
u\in L^\infty(0,\infty; H^1(\mathbb{T}^2)),\quad (\nabla(\nabla^\bot\cdot u), \nabla F,\sqrt\rho\dot u)\in L^2(0,\infty; L^2(\mathbb{T}^2)),\\[1mm] 
\sqrt t\sqrt\rho \dot u\in L^\infty_{\rm loc}(0,\infty; L^2(\mathbb{T}^2)),\quad \sqrt t (\nabla^\bot\cdot \dot u, \dot V)\in L^2_{\rm loc}(0,\infty; L^2(\mathbb{T}^2)),\\[1mm]
(\nabla^\bot\cdot u,\divg u)\in L^r_{\rm loc}(0,\infty; L^\infty(\mathbb{T}^2)),
\end{cases}
\end{equation}
where $1\leq p<\infty$, $1\leq r<2$ and $V=\divg u$. Moreover, if $P(\rho)=A\rho$ with $A>0$, then there exists at most one solution $(\rho,u)$, satisfying the regularity given in \eqref{1.13}, to the problem \eqref{1.1}--\eqref{1.4} on $\mathbb{T}^2\times[0,T]$.
 \end{thm}
 
 \begin{re}\label{re1.1} Indeed, the proof of  \cite[Proposition 3.3]{Danchin2023} gives a $\nu$-dependent estimate of the quantity on the left-hand side of \eqref{1.5}. Thus, for fixed and large $\nu\gg1$, the global existence theorem and the uniqueness result in \cite[Theorems 2.1, 2.2]{Danchin2023} remain valid. However, these $\nu$-independent estimates are technically needed for the convergence rate of the incompressible limit as $\nu\to\infty$ (see Theorem \ref{thm1.2} below).
\end{re}
 
Analogously to the proof in \cite{Danchin2023}, by virtue of the $\nu$-independent estimates established in Proposition \ref{pro1.2} and \cite[Propositions 3.1, 3.2]{Danchin2023} we can justify the incompressible limit (without convergence rate) from the compressible Navier-Stokes equations (\ref{1.1}) to the inhomogeneous incompressible Navier-Stokes equations, based on the classical compactness arguments. Formally,  we obtain by letting $\nu\to\infty$ that
\begin{align}
\begin{cases}\label{1.14}
\partial_t\eta+\divg(\eta v)=0 ,\\
\partial_t(\eta v)+\divg(\eta v\otimes v)-\mu\Delta v+\nabla\Pi = 0 ,\\
\divg v = 0,
\end{cases}
\end{align}
where $(\eta,v)$ is the limit of $(\rho,u)$ in some sense as $\nu\to\infty$. System (\ref{1.14}) describes the two-dimensional motion of an inhomogeneous incompressible viscous fluid (cf. \cite{AKM1990,Lions1996}).

When the initial data is more regular and the density is strictly away vacuum, Danchin-Mucha \cite{Danchin2019} proved the convergence rate of the incompressible limit from (\ref{1.1}) to (\ref{1.14}). In \cite{DFP2020}, the authors showed that if the initial density $\rho_0=1$, then the solution of the compressible Navier-Stokes equations (\ref{1.1}) converge to that of  the classical incompressible Navier-Stokes equations (i.e. (\ref{1.14}) with $\eta=1$), and the convergence rate of the incompressible limit was also obtained. However, to the authors' knowledge, the convergence rate of such incompressible limit remains unknown in the case that the solutions may contain vacuum states. A major difficulty  lies in  the lack of  the (weighted) bounds on  $u$ and $\dot u$, when the density vanishes and the vacuum appears. So, the second and main purpose of this paper is to derive the convergence rate of the incompressible limit as $\nu\to\infty$ in the presence of vacuum.

\begin{thm}\label{thm1.2} In addition to the conditions of Theorem \ref{thm1.1}, assume further that $\divg u_0=0$ and $\nabla\rho_0\in L^q$ with  some $q\in (2, \infty)$. Let $(\rho,u)$ and $(\eta,v)$ be the solutions of \eqref{1.1} and \eqref{1.14} with the same initial data $(\rho_0,u_0)$ on $\mathbb{T}^2\times[0,T]$, respectively. Then, there exists a positive constant $C(T)$, independent of $\nu=2\mu+\lambda$, such that for any $\nu\geq\nu_0$,
\begin{equation}
\begin{aligned}
&\sup_{0\leq t\leq T}\left(\|\sqrt\eta(\mathcal{P}u-v)(t)\|_{L^2}^2+\|(\rho-\eta)(t)\|_{L^2}^2\right)\\
&\quad+\int_0^T\left( \| \mathcal{P}u-v \|_{L^2}^2+\|\nabla ( \mathcal{P}u-v
)\|_{L^2}^2\right)\md t\leq C(T)\nu^{-\frac{1}{2}}
\end{aligned}
\end{equation}
and
\begin{equation}
\sup_{0\leq t\leq T}\|\nabla \mathcal{Q}u(t)\|_{L^2}^2+\nu\int_0^T\|\nabla\mathcal{Q} u\|_{H^1}^2\leq C(T)\nu^{-1},
\end{equation}
where  $\mathcal{P}:=\mathbb{Id}+\nabla(-\Delta)^{-1}\divg$ and $\mathcal{Q}:=-\nabla(-\Delta)^{-1}\divg$ are the Helmholtz projectors on divergence-free and potential vector fields, respectively.
 \end{thm}

\begin{re}\label{re1.2} If the initial density satisfies $\nabla\rho_0\in L^q$ with $2<q<\infty$, then the uniqueness of the solutions of the problem \eqref{1.1}--\eqref{1.4} can be shown to hold for the general pressure law $P(\rho)=A\rho^\gamma$ with $A>0$ and $\gamma\geq1$ in an easier manner than that in \cite{DFP2020,Danchin2023} by using the Euler coordinates (see, for example, \cite{Hu2021,Ger2011}). Similarly, the uniqueness of the solutions of system \eqref{1.14} can also be obtained in  an easier manner than that in \cite{Danchin2019-1}. 
\end{re}

We now make some comments on the proofs of Proposition \ref{pro1.2} and Theorem \ref{thm1.2}. We begin with the global $H^1$-estimate of the velocity and the uniform upper bound of the density stated in \cite[Propositions 3.1 and 3.2]{Danchin2023}, which are  actually global-in-time and $\nu$-independent (see Lemma  \ref{lem3.1}). As aforementioned, instead of (\ref{1.5}), we aim to prove (\ref{1.11}). The difference   between these two quantities  lies in the discrepancy between $\divg\dot u$ and $\dot V$ with $V=\divg u$. The advantage of handling $\dot V$ is that it enables us to  employ the ``effective viscous flux" $F$ defined in (\ref{1.10}) to deal with the terms associated with $\nu\nabla \divg u$ and $\nabla P$ (see (\ref{3.13}).  As observed by Hoff \cite{Ho1995}, the ``effective viscous flux" $F$ has better regularity than $\nu\divg u$ at this stage, and it holds that $\|\nabla F\|_{L^2(0,T;L^2)}$ (rather than $\nu\|\nabla \divg u\|_{L^2(0,T;L^2)}$) is uniformly bounded in $\nu$  (see (\ref{3.2})). Note that $\|\dot V\|_{L^2}-\|\nabla u\|_{L^4}^2 \leq  \|\divg \dot u\|_{L^2}\leq \|\dot V\|_{L^2}+\|\nabla u\|_{L^4}^2$ (see (\ref{3.14})). Unlike that in \cite{Danchin2023} where a Desjardins' inequality  (cf. \cite{Des1997}) was used,  to overcome the difficulty induced by the presence of vacuum,  we have to make use of a Poincar${\rm\acute{e}}$'s type inequality (cf. Lemma \ref{lem2.2}) to estimate the $L^p$-norm of the velocity. The technique of the ``curl-div" decomposition (cf. Lemma \ref{lem2.3}) also plays a key role in the entire analysis, since the divergence-free part and the potential part of the velocity exhibit different asymptotic behaviors with respect to large $\nu$. 

However, due to the lack of  the uniform-in-$\nu$  bound of  $\nu^{\frac{1}{4}}\|\nabla u\|_{L^4(0,T; L^4)}$, the last term  $N_6=-\nu\langle\dot V, \partial_iu\cdot\nabla u_i\rangle$ on the right-hand side of (\ref{3.13}) cannot be estimated directly. Indeed, this is the major difficulty lying in the proof of (\ref{1.11}). To this end,  noting that $\nu\dot V=(\nu\divg u)_t+u\cdot\nabla(\nu\divg u)$ and using the fact that  $\nu\divg u=F+P-\overline P$,  we find
\begin{align*}
N_6&=-\int \left[(\nu\divg u)_t+u\cdot\nabla(\nu\divg u)\right]\partial_i u\cdot\nabla u_i\mdx \\
&=-\int  F_t  \partial_i u\cdot\nabla u_i\mdx-\int \left(u\cdot\nabla F -\gamma P\divg u-\overline P_t\right)\partial_i u\cdot\nabla u_i\mdx,
\end{align*}
where the second term on the right-hand side can be  bounded by using the estimates stated in Lemma \ref{lem3.1}. For the first term, we integrate by parts to get
$$
\int  F_t   \partial_i u\cdot\nabla u_i\mdx
=\frac{\md}{\md t}\int  F   \partial_i u\cdot\nabla u_i\mdx+2\int   \left(F \partial_t u\cdot\nabla\divg u+\partial_t u\cdot\nabla u\cdot\nabla F \right)\mdx,
$$
By redirecting $u_t$ to the material derivative $\dot u$ and integrating by parts, we can use the Poincar${\rm\acute{e}}$'s type inequality (cf. Lemma \ref{lem2.2}) and the estimates achieved (cf. Lemma \ref{lem3.1}) to bound the second term. A key observation lying in  the first term is that 
$$
\sum_{i=1}^2 \partial_i u\cdot\nabla u_i=(\divg u)^2+2\nabla^\bot u_2\cdot\nabla u_1
$$
and
$$
\divg (\nabla^\bot u_2)=0,\quad \nabla^\bot\cdot(\nabla u_1)=0,
$$
The term $\langle F, (\divg u)^2\rangle$ will be bounded by using the ``effective viscous flux" $F$ again. More importantly,  based upon \cite[Theorem II.1]{Coifman1993} (also cf. Lemma \ref{lem2.5}) and the fact that BMO space and Hardy space are dual spaces (see  \cite{Meyer1975}), we have
$$\left|\langle F,\nabla^\bot u_2\cdot\nabla u_1\rangle\right|\leq \|\nabla F\|_{L^2}\|\nabla u\|_{L^2}^2\leq C\|\sqrt\rho \dot u\|_{L^2}\|\nabla u\|_{L^2}^2.
$$
With the above estimates at hand, we can obtain the desired estimates stated in Lemma \ref{lem3.2}. As a result, the bounds of $\|\rho\dot u\|_{L^r(0,T;L^q)}$ and $\|(\nabla^\bot\cdot u,\divg u)\|_{L^r(0,T;L^\infty)}$ can be deduced from the interpolation arguments. For more details, we refer to Section \ref{sec3}.

Based upon the uniform-in-$\nu$ estimates in Section \ref{sec3}, one can justify the incompressible limit as $\nu\to\infty$ in a similar manner as that in \cite{Danchin2023} by using the standard compactness arguments developed in \cite{Lions,Feireisl}. However, the convergence rate of the incompressible limit as $\nu\to\infty$ is a more subtle issue. To do this, we need the bound of $\|\nabla\rho\|_{L^q}$ with $2<q<\infty$, which strongly relies on the estimate of $\|\nabla u\|_{L^\infty}$ and will be achieved by solving a logarithmic inequality (see (\ref{4.5})) in an analogous manner as that in \cite{HLX2012}. We mention here that the bound of $\|\nabla u\|_{L^\infty}$ will be built upon the Beale-Kato-Majda's type inequality (cf. \cite{Kato,Huang}). 

The  fundamental idea of the proof of the convergence rate is to compute the discrepancy between the compressible and incompressible solutions by the (weighted) $L^2$-method. In fact, for the convergence of the velocities, we only need to compare  the divergence-free part of compressible velocity (i.e. $\mathcal{P}u$) with the incompressible velocity, since the potential part of compressible velocity (i.e. $\mathcal{Q} u$) decays sufficiently fast to zero as $\nu\to\infty$, due to Lemmas \ref{lem3.1}, \ref{lem3.3} and \ref{lem4.1}. However, some serious difficulties arise from the lack of the bounds on $(u,v)$ and $(u_t,v_t)$, due to the presence of vacuum.  To overcome these difficulties,   we first decompose the difference $\rho-\eta$ into two parts: 
$$
\rho-\eta=(\rho-\tilde\rho)+(\tilde\rho-\eta):=\varphi+\phi,
$$
where $\tilde\rho=\tilde\rho(x,t)$ is the transport of $\rho_0$ driven by the flow of the divergence-free part $\mathcal{P}u$ of the compressible velocity. Next, we need to carry out some $t$-growth and singular $t$-weighted estimates of $(\varphi,\phi)$ in $L^2$-norm. It is worth pointing out that the $t$-growth estimates of  $(\varphi,\phi)$ in $L^2$-norm ensure that the singular $t$-weighted estimates of $\|(\varphi,\phi)\|_{L^2}$ are meaningful. Moreover, due to the lower regularity of the initial velocity, by Lemmas \ref{lem2.2}, \ref{lem2.6} and \ref{lem3.2}   we only have the weighted estimates $\sqrt t\|(u_t,v_t)\|_{L^p}\in L^2(0,T)$, and thus, the singular $t$-weighted estimates of $\|(\varphi,\phi)\|_{L^2}$ are also technically needed  to  circumvent the difficulty caused by the lack of the bound on $\|u_t\|_{L^p}$ with $2<p<\infty$ (see, e.g. (\ref{4.19}) and (\ref{4.20})). For the same reasons, the treatment of the last term $M_8=-\langle \rho (\mathcal{Q}u)_t, \mathcal{P}u-v\rangle$ on the right-hand side of (\ref{4.18}) needs some mathematical tricks and will be done by making full use of the decay of $\mathcal{Q}u$. More precisely, noting that $\rho=\varphi+\phi+\eta$, we have
$$
M_8=-\int \left(\phi+\varphi\right) (\mathcal{Q}u)_t  \cdot (\mathcal{P}u-v)\mdx  -\int \eta (\mathcal{Q}u)_t \cdot (\mathcal{P}u-v)\mdx.
$$
The first term on the right-hand side can be bounded by using the singular $t$-weighted estimate of $\|(\varphi,\phi)\|_{L^2}$. To deal with the second term, we further write it in the form:
\begin{align*}
-\int \eta (\mathcal{Q}u)_t \cdot (\mathcal{P}u-v)\mdx&=-\frac{\md}{\md t}\int \eta (\mathcal{Q}u) \cdot (\mathcal{P}u-v)\mdx-\int v\cdot\nabla\eta (\mathcal{Q}u) \cdot (\mathcal{P}u-v)\mdx\\
&\quad +\int \eta  (\mathcal{Q}u)  \cdot (\mathcal{P}u_t-v_t)\mdx, 
\end{align*}
where the second term on the right-hand side can be well controlled by  the fast decay of $\mathcal{Q}u$. For the third term, noting that $\mathcal{P}u_t=u_t-\mathcal{Q}u_t$ and $\eta=\eta^{\frac{1}{2}}(\rho-\varphi-\phi)^{\frac{1}{2}}$, we have
\begin{align*}
\int \eta  (\mathcal{Q}u)  \cdot (\mathcal{P}u_t-v_t)\mdx&=-\int \eta  (\mathcal{Q}u)  \mathcal{Q}u_t\mdx+\int \eta  (\mathcal{Q}u)  (u_t-v_t)\mdx\\
&\leq -\int \eta  (\mathcal{Q}u)  \mathcal{Q}u_t\mdx + \int \eta^{\frac{1}{2}}
\rho^{\frac{1}{2}} |\mathcal{Q}u| | u_t-v_t|\mdx\\
&\quad 
+ \int \eta^{\frac{1}{2}}
(|\phi|+|\varphi|)^{\frac{1}{2}} |\mathcal{Q}u| | u_t-v_t|\mdx.
\end{align*}
The above decomposition, together with the estimates of $\|(\sqrt\rho u_t,\sqrt\eta v_t)\|_{L^2}$ and the singular $t$-weighted estimates of $\|(\varphi,\phi)\|_{L^2}$,  enables us to bound the right-hand side  by the fast decay of $\mathcal{Q}u$. For the detailed analysis, we refer to Section \ref{sec4}.

The rest  of this paper is organized as follows.  In Section \ref{sec2}, we collect
some useful inequalities and known facts, which will be  used in later analysis.  Section \ref{sec3} is devoted to the global $\nu$-independent and $t$-weighted estimates stated in Proposition \ref{pro1.2} and the proof of Theorem \ref{thm1.1}. The convergence rate of the incompressible limit as $\nu\to\infty$ in Theorem \ref{thm1.2} will be shown in Section \ref{sec4}.

\section{Preliminaries}\label{sec2}
In this section, we recall some known results and elementary inequalities.  Without loss of generality, we assume that  $|\mathbb{T}^2|=1$. For any $1\leq p\leq \infty$, we simply denote by 
$$
L^p:=L^p(\mathbb{T}^2),\ \ H^1:=H^1(\mathbb{T}^2),\ \ W^{1,p}:=W^{1,p}(\mathbb{T}^2)\ \
{\text{and}}\ \ 
\int f(x)\mdx:=\int_{\mathbb{T}^2} f(x)\mdx.
$$

We begin with the well-known Gagliardo-Nirenberg-Sobolev's inequalities in 2D (see  \cite{Nirenberg}).
\begin{lem}\label{lem2.1} For  $p\geq2$, $q>1$ and $r>2$, assume that $f\in H^1 $ and $g\in L^q \cap W^{1,r} $. Then there exist positive constants $C$, $C_1$ and $C_2$, depending on $p$, $q$ and $r$, such that, 
\begin{align}\label{2.1}
\|f\|_{L^p}\leq C \|f\|_{L^2}^{\frac{2}{p}}\|\nabla f\|_{L^2}^{\frac{p-2}{p}}+C_1\|f\|_{L^2},
\end{align}
and
\begin{align}\label{2.2}
\|g\|_{L^\infty}\leq C\|g\|_{L^q}^{\frac{q(r-2)}{2r+q(r-2)}}\|\nabla g\|_{L^r}^{\frac{2r}{2r+q(r-2)}}+C_2\|g\|_{L^2}.
\end{align}
Moreover, if $\bar{f}=0$ (resp. $\bar{g}=0$), then $C_1=0$ (resp. $C_2=0$).
\end{lem}

In order to estimate the $L^p$-norm of the velocity, we need to use the Poincar\'{e}-type inequality, which can be found in  \cite[Lemma 3.2]{Feireisl}.
\begin{lem}\label{lem2.2}Assume that there are two positive constants $M_1$ and $M_2$ such that the non-negative function $\rho=\rho(x)$ satisfies
\begin{equation}\label{2.3}
0<M_{1} \leq \int \rho(x) \mdx ,\quad  \int \rho^{\theta}(x) \mdx \leq M_{2}\quad{\text with }\quad \theta>1,
\end{equation}
Then for any $p\geq2$, there exists a positive constant $C$, depending only on $M_{1}, M_{2}$, $\theta$ and $p$, such that for any $v \in H^{1} $,
\begin{align}\label{2.4}
\|v\|_{L^{p}}^{2} \leq C\left(\|\sqrt \rho v\|^{2}_{L^2}+\|\nabla v\|_{L^{2}}^{2}\right).
\end{align}
\end{lem}

Obviously, the condition (\ref{2.3}) is automatically satisfied for the general pressure law $P(\rho)=A\rho^\gamma$ with $\gamma\geq1$, provided the density is uniformly bounded on $\mathbb{T}^2$. Next, we recall the standard ``div-curl" estimate (see \cite{Brezis1974}), which implies that the norm of $\|\nabla u\|_{L^p}$ with $1<p<\infty$ can be well controlled by the $L^p$-norms of   
of ${\rm div} u$ and $\nabla^\bot\cdot u$.
\begin{lem}\label{lem2.3}
Let $k\geq 0$ be an integer and  $1<p<\infty$, assume that $u\in W^{k+1,p} $. Then there exists a positive constant $C=C(k,p)$ such that
\begin{align}\label{2.5}
\|\nabla u\|_{W^{k,p}}\leq C(\|\divg u\|_{W^{k,p}}+\|\nabla^\bot\cdot u\|_{W^{k,p}}).
\end{align}
\end{lem}

It is  known that the $W^{k,p}$-estimates stated in  (\ref{2.5}) is not valid for the endpoint case $p=\infty$. So, to deal with $\|\nabla u\|_{L^\infty}$, we need to use the
  Beale-Kato-Majda's type inequality (see \cite{Kato,Huang}).
\begin{lem}\label{lem2.4}
For $2<q<\infty$, assume that $\nabla u\in  W^{1,q} $. Then there exists a positive constant $C=C(q)$ such that
\begin{align}\label{2.6}
\|\nabla u\|_{L^\infty}\leq C(\|\divg u\|_{L^\infty}+\|\nabla^\bot\cdot u\|_{L^\infty})\ln(e+\|\nabla^2 u\|_{L^q})+C\|\nabla u\|_{L^2}+C.
\end{align}
\end{lem}

Let $(-\Delta)^{-1}$ be the inverse of Laplacian operator with zero mean on $\mathbb{T}^2$ and $\mathcal{R}_i:=(-\Delta)^{-1/2}\partial_i$ be the usual Riersz transform on $\mathbb{T}^2$.  Then, the Hardy space $\mathcal{H}^1(\mathbb{T}^2)$ and the BMO space ${\text{BMO}}(\mathbb{T}^2)$ are defined respectively by (see  \cite{Meyer1975}):
\begin{align*}
\mathcal{H}^1:=\left\{f\in L^1(\mathbb{T}^2):\|f\|_{\mathcal{H}^1}:=\|f\|_{L^1}+\|\mathcal{R}_1f\|_{L^1}+\|\mathcal{R}_2f\|_{L^1}<\infty,~\bar{f}=0\right\},
\end{align*}
and
\begin{align*}
\text{BMO}:=\left\{f\in L_{\rm loc}^1(\mathbb{T}^2):\|f\|_{\text{BMO}}<\infty\right\},
\end{align*}
with
\begin{align*}
\|f\|_{\text{BMO}}:=\sup_{x\in\mathbb{T}^2,r\in(0,1)}\frac{1}{|\Omega_r(x)|}\int_{\Omega_r(x)}\left|f(y)-\frac{1}{|\Omega_r(x)|}\int_{\Omega_r(x)}f(z)\mathrm{d}z\right|\md y,
\end{align*}
where  $\Omega_r(x)=\mathbb{T}^2\cap B_r(x)$ and $B_r(x)$ is a ball with center $x$ and radius
$r$. It follows from \cite{Meyer1975} that $\mathcal{H}^1(\mathbb{T}^2)$ and ${\text{BMO}}(\mathbb{T}^2)$ are dual spaces. In view of  \cite[Theorem II.1]{Coifman1993}, we also have

\begin{lem}\label{lem2.5} Assume that $(u,v)\in L^2$ satisfies  $\nabla^\bot \cdot u=0$ and $\divg v=0$ in the sense of distribution. Then, it holds that $u\cdot v \in \mathcal{H}^1(\mathbb{T}^2)$  and 
\begin{equation}\label{2.7}
 \|u\cdot v\|_{\mathcal{H}^1}\leq \|u\|_{L^2}\|v\|_{L^2}.
\end{equation}
\end{lem}

Finally, we end this section with the existence theorem for the inhomogeneous incompressible Navier-Stokese equations (\ref{1.14}). In fact, the global well-posedness of strong solutions of system (\ref{1.14}) was shown in \cite{Danchin2019-1} under the weaker conditions that $\rho_0\in L^\infty$ and $u_0\in H^1$. If the initial density  additionally fulfills $\nabla\rho_0\in L^q $ with $2<q<\infty$, then it is easily deduced that $\nabla\eta\in L^\infty_{\rm loc}(0,\infty; L^q )$ in a similar manner as that in Lemma \ref{lem3.3}. It is worth mentioning here that the $W^{1,q}$-estimate of the density plays an important role in the derivation of the convergence rate of the incompressible limit.. 

\begin{lem}\label{lem2.6}
Assume that  $(\rho_0,u_0)$ satisfies the conditions of Theorem \ref{thm1.2}. Then for any $T>0$, there exists a global unique strong solution $(\eta,v)$ to the problem \eqref{1.14}
 with the initial data $(\rho_0,u_0)$ on $\mathbb{T}^2\times[0,T]$, such that for $p\in[1,\infty)$ and $q\in[1,2)$,
 \begin{equation}\label{2.8}
 \begin{cases}
0\leq \eta\in L^\infty(0,\infty; W^{1,q}(\mathbb{T}^2))\cap C([0,\infty); L^p(\mathbb{T}^2)),\\[1mm]
 \sqrt\eta v\in C([0,\infty); L^2(\mathbb{T}^2)),\quad v\in L^\infty(0,\infty; H^1(\mathbb{T}^2)), \\[1mm]
(\sqrt\eta v_t,\nabla^2v )\in L^2(0,\infty; L^2(\mathbb{T}^2)),\quad \nabla v\in L^1_{\rm loc}(0,\infty; L^\infty(\mathbb{T}^2)),\\[1mm]
\sqrt t(\sqrt\eta v_t, \nabla^2v)\in  L_{\rm loc}^\infty (0,\infty; L^2(\mathbb{T}^2)),\quad  \sqrt t\nabla v_t\in L_{\rm loc}^2 (0,\infty; L^2(\mathbb{T}^2)).
\end{cases} 
\end{equation}
 \end{lem}

\section{Global $t$-weighted estimates}\label{sec3}
As aforementioned, the proof of Proposition \ref{pro1.2} is built upon the global uniform-in-$\nu$ estimates in \cite[Propostions 3.1 and 3.2]{Danchin2023}. For the sake of convenience, we summarize   these estimates in the following lemma. For simplicity, we assume that $\nu\geq1$ throughout of this paper.

\begin{lem}\label{lem3.1}
Let the conditions of Proposition \ref{pro1.2} be in force. For any $0<T<\infty$, assume that $(\rho,u)$ is a smooth solution of the problem \eqref{1.1}--\eqref{1.4} on $\mathbb{T}^2\times[0,T]$. There exists a positive number $\nu_0>0$, depending only on $A$, $\gamma$, $\mu$,  $\|\rho_0\|_{L^\infty}$, $\|u_0\|_{H^1}$ and $K$, such that if $\nu\geq\nu_0$, then
\begin{equation}\label{3.1}
0\leq\rho (x,t)\leq C,\quad\forall\ x\in\mathbb{T}^2,\ t\in[0,T],
\end{equation}
and
\begin{equation}\label{3.2}
\begin{aligned}
&\sup_{0\leq t\leq T}\left(\|\sqrt\rho u\|_{L^2}^2+ \|\nabla^\bot\cdot u\|_{L^{2}}^{2} +\nu\|\divg u\|_{L^2}^2\right)(t)+ \int_0^T\left(\|\nabla u\|_{L^2}^2+\|\sqrt\rho\dot u\|_{L^2}^2\right)\md t\\
&\qquad+ \int_0^T\left(\|\nabla(\nabla^\bot\cdot u)\|_{L^2}^2+\|\nabla F\|_{L^2}^2+  \frac{1}{\nu}\|P-\bar P\|_{L^2}^2+\|\nabla u\|_{L^4}^4 \right)\md t\leq C,
\end{aligned}
\end{equation}
where $C>0$ is a   generic positive constant independent of $\nu$ and $T$.
\end{lem}
\pf Indeed, noting that $\divg u=\nu^{-1}(F+P-\overline P)$ and  $\nabla^\bot\cdot\mathcal{Q}u=0$, we have
$$
\nu\|\divg u\|_{L^2}^2\leq \nu^{-1}\left(\|F\|_{L^2}^2+\|P-\overline P\|_{L^2}^2\right)
$$
and
$$
\|\nabla^k(\nabla^\bot\cdot u)\|_{L^2}=\|\nabla^k(\nabla^\bot\cdot\mathcal{P}u)\|_{L^2},\quad k=0,1.
$$
So, combining  \cite[Propostions 3.1 and 3.2]{Danchin2023} with  \cite[ (3.47)]{Danchin2023} proves (\ref{3.1}) and (\ref{3.2}).\hfill$\square$

\vskip 2mm

The next lemma is concerned with the globally $\nu$-independent and $t$-weighted estimates of the solutions. To do this, from now on, we assume that $\nu\geq\nu_0$. Moreover, for simplicity, we denote by $C$ and $C_i$ ($i=1,2,\ldots$) the positive constants which may depend on $T$, $A$, $\gamma$, $\mu$, $K$ and the initial norms, but not on $\nu$. 

\begin{lem}\label{lem3.2}
Let the conditions of  Proposition \ref{pro1.2}  be satisfied. Then,  
\begin{equation}
\label{3.3}
\sup_{0\leq t\leq T} \left(t \|\sqrt\rho \dot{u}(t)\|_{L^2}^2 \right) +\int_0^Tt\left(\mu\|\nabla^\bot\cdot \dot u\|_{L^2}^2+\nu\| \dot V\|_{L^2}^2\right)\md t 
\leq C(T),
\end{equation}
where $V:=\divg u$ and $\dot V:=\divg u_t+u\cdot\nabla \divg u$.
\end{lem}
\pf  To clarify the proof, we divide it into several steps.

\vskip 2mm

\underline{\bf Step 1. Differential identity for $\|\sqrt\rho \dot{u}\|_{L^2}^2$}

\vskip 2mm

To begin, we first rewrite (\ref{1.1})$_2$ in the form:
 \begin{equation}\label{3.4}
 \rho\dot u+\nabla P(\rho )=\mu\nabla^\bot(\nabla^\bot\cdot u)+
\nu\nabla\divg u\quad{\text{with}}\quad \nu=2\mu+\lambda.
\end{equation}
Then, operating $\partial_t+\divg(u\cdot )$ to both sides of (\ref{3.4}) and  multiplying it by $\dot u$ in $L^2$, we obtain 
\begin{equation}
\begin{aligned}\label{3.5}
\frac{1}{2}\frac{\md}{\md t} \|\sqrt\rho \dot{u}\|_{L^2}^2&=\mu \int\dot{u}\cdot\left(\nabla^\bot(\nabla^\bot\cdot u_t)
+\divg\left(u\left(\nabla^\bot(\nabla^\bot\cdot u)\right)\right)\right)\mdx\\
&\quad+\nu\int\dot{u}\cdot\left(\nabla \divg u_t+\divg 
(u\nabla\divg  u)\right)\mdx\\
&\quad -\int\dot{u}\cdot\left(\nabla P_t+\divg 
(u\nabla P)\right)\mdx:=I_1+I_2+I_3.
\end{aligned}
\end{equation}

By direct calculations, we have
\begin{equation}\label{3.6}
\nabla^\bot(\nabla^\bot\cdot u_t)
+\divg\left(u\left(\nabla^\bot(\nabla^\bot\cdot u)\right)\right)=\nabla^\bot(\nabla^\bot\cdot \dot u)-\divg \left((\nabla^\bot\cdot u)(\nabla^\bot \otimes u)\right),
\end{equation}
where $\nabla^\bot\otimes u$ is the tensor product of $\nabla^\bot$ and $u$ defined by
$$
\nabla^\bot\otimes u=(-\partial_2 u, \partial_1 u),\quad \divg (\nabla^\bot\otimes u)=(-\partial_2\divg u,\partial_1\divg u),
$$
and
$$\divg \left((\nabla^\bot\cdot u)(\nabla^\bot \otimes u)\right)= \left(-\divg  ((\nabla^\bot \cdot u)\partial_2 u ), \divg  ((\nabla^\bot\cdot u)\partial_1u )\right).$$
Let  $A:B:=a_{ij}b_{ij}$ for $A=(a_{ij})_{2\times 2}$ and $B=(b_{ij})_{2\times2}$. Then, it follows from (\ref{3.6}) that
\begin{equation}\label{3.7}
I_1= -\mu \|\nabla^\bot\cdot \dot u\|_{L^2}^2+\mu\int (\nabla^\bot \cdot u) \nabla\dot u:(\nabla^\bot\otimes u)\mdx.
\end{equation}

Recalling that $V=\divg u$ and $\dot V=V_t+u\cdot\nabla V$, 
we find
\begin{equation}\label{3.8}
\begin{aligned}
\nabla \divg u_t+\divg 
(u\nabla\divg  u)&=\nabla (V_t+u\cdot\nabla V)-\nabla (u\cdot\nabla \divg u)+\divg (u\nabla \divg u)\\
&=\nabla\dot V-\nabla u\cdot\nabla \divg u+\divg u\nabla \divg u\\
&=\nabla\dot V-\partial_i ((\divg u) \nabla u_i)+2(\divg u)\nabla \divg u,
\end{aligned}
\end{equation}
where we have used Einstein convention that  repeated indices denote summation. 
Due to (\ref{1.1})$_1$, it is easily seen that
\begin{equation}\label{3.9}
\begin{aligned}
\nabla P_t+\divg (u\nabla P)&=-\nabla (\gamma P\divg u)-\nabla  u\cdot\nabla P+(\divg u)\nabla P\\
&=-\nabla (\gamma P\divg u)-\partial_i(P\nabla  u_i)+(\divg u)\nabla P+P\nabla \divg u.
\end{aligned}
\end{equation}

Based upon (\ref{3.8}) and (\ref{3.9}), we infer from integration by parts that
\begin{equation}\label{3.10}
\begin{aligned}
I_2+I_3&=-\int \divg \dot u\left(\nu \dot V+\gamma P\divg u\right)\md x+ \int  \partial_i \dot u \cdot  ((\nu\divg u-P) \nabla u_i)\md x\\
&\quad +\int (\divg u) \dot u\cdot \nabla (\nu\divg u- P) \md x+\int  (\nu\divg u-P) \dot u\cdot \nabla \divg u \md x.
\end{aligned}
\end{equation}

Noting that
\begin{equation}
\divg \dot u=(\divg u)_t+\divg (u\cdot\nabla u)=\dot V-u\cdot\nabla\divg u+\divg (u\cdot\nabla u)=\dot V+\partial_i u\cdot\nabla u_i, \label{3.11}
\end{equation}
which, inserted into the first term on the right-hand side of (\ref{3.10}), gives
\begin{equation}\label{3.12}
\begin{aligned}
I_2+I_3&=-\nu\int |\dot V|^2\md x - \nu \int\dot V \partial_iu\cdot\nabla u_i\md x-\gamma \int P (\divg \dot u)(\divg u) \md x\\
&\quad+ \int  \partial_i \dot u \cdot  ((\nu\divg u-P) \nabla u_i)\md x+\int (\divg u) \dot u\cdot \nabla (\nu\divg u- P) \md x\\
&\quad +\int  (\nu\divg u-P) \dot u\cdot \nabla \divg u \md x.
\end{aligned}
\end{equation}

Thus, plugging (\ref{3.7}) and (\ref{3.12}) into (\ref{3.5}), we arrive at
\begin{equation}
\begin{aligned}\label{3.13}
&\frac{1}{2}\frac{\md}{\md t} \|\sqrt\rho \dot{u}\|_{L^2}^2 +\left(\mu\|\nabla^\bot\cdot \dot u\|_{L^2}^2+\nu\| \dot V\|_{L^2}^2\right)\\
&\quad=\mu\int (\nabla^\bot \cdot u) \nabla\dot u:(\nabla^\bot\otimes u)\mdx  
- \gamma\int P (\divg \dot u)(\divg u) \md x\\
&\qquad + \int  \partial_i \dot u \cdot  ((\nu\divg u-P) \nabla u_i)\md x  +\int (\divg u) \dot u\cdot \nabla (\nu\divg u- P) \md x\\
&\qquad +\int  (\nu\divg u-P) \dot u\cdot \nabla \divg u \md x - \nu \int\dot V \partial_iu\cdot\nabla u_i\md x  :=N_1+\ldots+N_6.
\end{aligned}
\end{equation}

\vskip 2mm

\underline{\bf Step 2. The estimates for the terms $N_1$--$N_5$}

\vskip 2mm

In order to deal with the right-hand side of (\ref{3.13}), we first infer from (\ref{3.11}) that
\begin{equation}\label{3.14}
\|\dot V\|_{L^2}-\|\nabla u\|_{L^4}^2 \leq  \|\divg \dot u\|_{L^2}\leq \|\dot V\|_{L^2}+\|\nabla u\|_{L^4}^2,
\end{equation}
so that, by (\ref{2.5}) we obtain
\begin{equation}\label{3.15}
\begin{aligned}
|N_1| &\leq C\|\nabla u\|_{L^4}^2\|\nabla \dot u\|_{L^2}\leq C\|\nabla u\|_{L^4}^2\left(\|\nabla^\bot\cdot \dot u\|_{L^2}+\|\divg \dot u\|_{L^2}\right)\\
&\leq  \frac{1}{16}\left(\mu\|\nabla^\bot\cdot\dot u\|_{L^2}^2+ \nu \|\dot V\|_{L^2}^2\right)+C \|\nabla u\|_{L^4}^4 .
\end{aligned}
\end{equation}

Owing to (\ref{3.1}), (\ref{3.2}) and (\ref{3.14}), it is easily seen that
\begin{equation}\label{3.16}
|N_2|\leq C\|\divg \dot u\|_{L^2}\|\divg  u\|_{L^2}\leq \frac{\nu}{16}\|\dot V\|_{L^2}^2+C\left(1+\|\nabla u\|_{L^4}^4\right).
\end{equation}

To deal with the third term $N_3$, we first observe from (\ref{3.4}) that
\begin{equation}\label{3.17}
\Delta F=\divg(\rho\dot u)\quad{\text{and}}\quad \Delta (\nabla^\bot \cdot u)=\nabla^\bot\cdot(\rho\dot u),
\end{equation}
which, together with (\ref{3.1}) and the Poincar${\rm\acute{e}}$ inequality $\|F\|_{L^2}\leq C\|\nabla F\|_{L^2}$ (due to the fact that the mean of $F$ is zero (i.e. $\overline F=0$)), gives
\begin{equation}\label{3.18}
\|F\|_{L^2}+\|\nabla F\|_{L^2}+\|\nabla(\nabla^\bot\cdot u)\|\leq C\|\sqrt\rho\dot u\|_{L^2}.
\end{equation}

Using (\ref{2.5}), (\ref{3.1}), (\ref{3.2}), (\ref{3.14}) and (\ref{3.18}), by (\ref{1.10}) we find
\begin{equation}\label{3.19}
\begin{aligned}
|N_3|&= \int  \partial_i \dot u \cdot  ((F-\overline P) \nabla u_i)\md x\leq C\|\nabla \dot u\|_{L^2}\left(\|F\|_{L^4}\|\nabla u\|_{L^4}+\|\nabla u\|_{L^2}\right)\\
&\leq C\left(\|\nabla^\bot\cdot \dot u\|_{L^2}+\|\divg \dot u\|_{L^2} \right)\left(1+\|F\|_{H^1}\|\nabla u\|_{L^4}\right)\\
&\leq C\left(\|\nabla^\bot\cdot \dot u\|_{L^2}+\|\dot V\|_{L^2}+\|\nabla u\|_{L^4}^2\right)\left(1+\|\sqrt\rho \dot u\|_{L^2}\|\nabla u\|_{L^4}\right)\\
&\leq \frac{1}{16}\left(\mu\|\nabla^\bot\cdot\dot u\|_{L^2}^2+\nu\|\dot V\|_{L^2}^2\right)+C\left(1+\|\nabla u\|_{L^4}^4+\|\sqrt\rho \dot u\|_{L^2}^4\right).
\end{aligned}
\end{equation}

In a similar manner, by virtue of (\ref{2.4}) we deduce
\begin{equation}\label{3.20}
\begin{aligned}
|N_4|&\leq C\|\divg u\|_{L^4}\|\dot u\|_{L^4}\|\nabla F\|_{L^2}\leq  C\|\divg u\|_{L^4}\left(\|\sqrt\rho \dot u\|_{L^2}+\|\nabla \dot u\|_{L^2}\right)\|\sqrt\rho \dot u\|_{L^2}\\
&\leq  C\|\divg u\|_{L^4}\left(\|\sqrt\rho \dot u\|_{L^2}+\|\nabla^\bot\cdot \dot u\|_{L^2}+\|\dot V\|_{L^2}+\|\nabla u\|_{L^4}^2\right)\|\sqrt\rho \dot u\|_{L^2}\\
&\leq \frac{1}{16}\left(\mu\|\nabla^\bot\cdot\dot u\|_{L^2}^2+\nu\|\dot V\|_{L^2}^2\right)+C\left(1+\|\nabla u\|_{L^4}^4+ \|\sqrt\rho \dot u\|_{L^2}^4\right).
\end{aligned}
\end{equation}

After integrating by parts, we infer from (\ref{2.4}), (\ref{2.5}), (\ref{3.1}), (\ref{3.2}), (\ref{3.14}) and (\ref{3.18}) that
\begin{equation}\label{3.21}
\begin{aligned}
|N_5|&=\left|\int (F-\overline P)\dot u\cdot\nabla\divg u\mdx\right| =\left|\int (\divg u)\left(\dot u\cdot\nabla F +F\divg \dot u- \overline P \divg\dot u  \right)\mdx\right|\\
&\leq C\|\divg u\|_{L^4}\|\dot u\|_{L^4}\|\nabla F\|_{L^2}+  C\left(\|\divg u\|_{L^4}\|F\|_{L^4}+\|\divg u\|_{L^2}\right)\|\divg \dot u\|_{L^2}\\
&\leq  C\|\divg u\|_{L^4}\left(\|\sqrt\rho \dot u\|_{L^2}+\|\nabla^\bot\cdot \dot u\|_{L^2}+\|\dot V\|_{L^2}+\|\nabla u\|_{L^4}^2\right)\|\sqrt\rho \dot u\|_{L^2}\\
&\quad +C\left(1+\|\divg u\|_{L^4}\|\sqrt\rho\dot u\|_{L^2}  \right)\left( \|\dot V\|_{L^2}+\|\nabla u\|_{L^4}^2\right)\\
&\leq \frac{1}{16}\left(\mu\|\nabla^\bot\cdot\dot u\|_{L^2}^2+\nu\|\dot V\|_{L^2}^2\right)+C\left(1+\|\nabla u\|_{L^4}^4+\|\sqrt\rho \dot u\|_{L^2}^4\right).
\end{aligned}
\end{equation}

\vskip 2mm

\underline{\bf Step 3. The estimate for the term $N_6$}

\vskip 2mm
We now deal with the most troublesome term $N_6$. First, it follows from (\ref{1.10}) that
\begin{equation}\label{3.22}
\begin{aligned}
N_6&=-\int \left[(\nu\divg u)_t+u\cdot\nabla(\nu\divg u)\right]\partial_i u\cdot\nabla u_i\mdx \\
&=-\int \left[(F+P-\overline P)_t+u\cdot\nabla(F+P-\overline P)\right]\partial_i u\cdot\nabla u_i\mdx\\
&=-\int  F_t  \partial_i u\cdot\nabla u_i\mdx-\int \left(u\cdot\nabla F\right)\partial_i u\cdot\nabla u_i\mdx\\
&\quad  -\int \left[(P-\overline P)_t+u\cdot\nabla(P-\overline P)\right]\partial_i u\cdot\nabla u_i\mdx:= N_6^1+N_{6}^2+N_{6}^3.
\end{aligned}
\end{equation}

To consider the first term on the right-hand side of (\ref{3.22}), we   rewrite it as follows, 
\begin{equation}\label{3.23}
\begin{aligned}
N_6^1&=-\frac{\md}{\md t}\int  F   \partial_i u\cdot\nabla u_i\mdx+\int  F  \left( \partial_i \partial_t u\cdot\nabla u_i+\partial_iu\cdot\nabla\partial_t u_{i}\right)\mdx\\
&=-\frac{\md}{\md t}\int  F   \partial_i u\cdot\nabla u_i\mdx-2\int   F \partial_t u\cdot\nabla\divg u\mdx-2\int \partial_t u\cdot\nabla u\cdot\nabla F \mdx.
\end{aligned}
\end{equation}

Keeping in mind that $u_t=\dot u-u\cdot\nabla u$, we have
\begin{align*}
&-2\int   F \partial_t u\cdot\nabla\divg u\mdx= -2\int   F \dot u\cdot\nabla\divg u\mdx+2\int   F u\cdot\nabla u\cdot\nabla\divg u\mdx\\
&\quad= 2\int  (\divg u) \left( \dot u\cdot\nabla F + F\divg \dot u \right) \mdx-2\int   (\divg u) u\cdot\nabla u\cdot\nabla F\mdx\\
&\qquad -2\int   F(\divg u) \partial_i u\cdot\nabla u_i\mdx+\int  (\divg u)^2  u\cdot\nabla  F\mdx+\int  F (\divg u)^3 \mdx
\end{align*}
and
$$
-2\int \partial_t u\cdot\nabla u\cdot\nabla F \mdx =-2\int \dot u\cdot\nabla u\cdot\nabla F \mdx +2\int  u\cdot\nabla u\cdot\nabla u\cdot\nabla F \mdx,
$$
which, inserted into (\ref{3.23}), gives rise to
\begin{equation}\label{3.24}
\begin{aligned}
N_{6}^1 
&=-\frac{\md}{\md t}\int  F   \partial_i u\cdot\nabla u_i\mdx+2\int  (\divg u) \left( \dot u\cdot\nabla F + F\divg \dot u \right) \mdx\\
&\quad -2\int   (\divg u) u\cdot\nabla u\cdot\nabla F\mdx 
 -2\int   F(\divg u) \partial_i u\cdot\nabla u_i\mdx\\
 &\quad +\int  (\divg u)^2  u\cdot\nabla  F\mdx+\int  F (\divg u)^3 \mdx\\
 &\quad -2\int \dot u\cdot\nabla u\cdot\nabla F \mdx +2\int  u\cdot\nabla u\cdot\nabla u\cdot\nabla F \mdx\\
 &:=-\frac{\md}{\md t}\int  F   \partial_i u\cdot\nabla u_i\mdx+ N_{6,1}^1+\ldots+N_{6,7}^1.
\end{aligned}
\end{equation}

Similarly to the derivations of (\ref{3.19}) and (\ref{3.21}), we find
\begin{equation}\label{3.25}
\begin{aligned}
N_{6,1}^1+N_{6,6}^1&\leq C\|\nabla u\|_{L^4}\|\dot u\|_{L^4}\|\nabla F\|_{L^2} +C\|\nabla u\|_{L^4}\|F\|_{L^4}\|\divg \dot u\|_{L^2}\\
&\leq \frac{1}{16}\left(\mu\|\nabla^\bot\cdot\dot u\|_{L^2}^2+\nu\|\dot V\|_{L^2}^2\right) +C\left(1+\|\nabla u\|_{L^4}^4+ \|\sqrt\rho \dot u\|_{L^2}^4\right).
\end{aligned}
\end{equation}

It is easily deduced from (\ref{2.1}) and (\ref{3.18}) that
\begin{equation}\label{3.26}
N_{6,3}^1+N_{6,5}^1\leq C\|F\|_{L^4}\|\nabla u\|_{L^4}^3\leq C\|\nabla F\|_{L^2}\|\nabla u\|_{L^4}^3\leq C\left( \|\nabla u\|_{L^4}^4+\|\sqrt\rho \dot u\|_{L^2}^4 \right).
\end{equation}

Based upon (\ref{2.2}), (\ref{2.4}), (\ref{3.2}) and (\ref{3.18}), 
the remaining terms can  be bounded by 
\begin{equation}\label{3.27}
\begin{aligned}
&N_{6,2}^1+N_{6,4}^1+N_{6,7}^1\leq C\|u\|_{L^\infty}\|\nabla u\|_{L^4}^2\|\nabla F\|_{L^2}\\
&\quad\leq C\left(\|\sqrt\rho u\|_{L^2}+\|\nabla u\|_{L^2}+\|\nabla u\|_{L^4}\right) \|\nabla u\|_{L^4}^2\|\sqrt\rho \dot u\|_{L^2}\\
&\quad\leq C\left(1+ \|\nabla u\|_{L^4}^4+ \|\sqrt\rho \dot u\|_{L^2}^4\right).
\end{aligned}
\end{equation}

Thus, substituting (\ref{3.25})--(\ref{3.27}) into (\ref{3.24}), we obtain
\begin{equation}\label{3.28}
\begin{aligned}
N_{6}^1 
 &\leq -\frac{\md}{\md t}\int  F   \partial_i u\cdot\nabla u_i\mdx+\frac{1}{16}\left(\mu\|\nabla^\bot\cdot\dot u\|_{L^2}^2+\nu\|\dot V\|_{L^2}^2\right)\\
  &\quad+ C\left(1+\|\nabla u\|_{L^4}^4+\|\sqrt\rho \dot u\|_{L^2}^4\right).
\end{aligned}
\end{equation}

Analogously to the derivation of (\ref{3.27}), we easily see that
\begin{equation}\label{3.29}
N_{6}^2 \leq C\left( 1+\|\nabla u\|_{L^4}^4+ \|\sqrt\rho \dot u\|_{L^2}^4\right).
\end{equation}

Due to (\ref{1.1})$_1$, one has $P_t=-u\cdot\nabla P-\gamma P\divg u$. So, it follows from (\ref{3.1}) and (\ref{3.2}) that
$$
\left|\overline P'(t)\right|\leq \left|\int P(\rho)_t \mdx\right| \leq C\int |P||\divg u|\mdx\leq C\|\nabla u\|_{L^2}\leq C, 
$$
which, together with (\ref{3.2}), yields
\begin{equation}\label{3.30}
N_{6}^3\leq  C\left(\|\nabla u\|_{L^2}^2+\|\nabla u\|_{L^2}\|\nabla u\|_{L^4}^2\right)\leq C\left( 1+  \|\nabla u\|_{L^4}^4\right).
\end{equation}

Hence, plugging (\ref{3.28})--(\ref{3.30}) into (\ref{3.22}), we obtain
\begin{equation}\label{3.31}
\begin{aligned}
N_{6} 
 &\leq -\frac{\md}{\md t}\int  F   \partial_i u\cdot\nabla u_i\mdx +\frac{1}{16}\left(\mu\|\nabla^\bot\cdot\dot u\|_{L^2}^2+\nu\|\dot V\|_{L^2}^2\right)\\
 &\quad + C\left( 1+ \|\nabla u\|_{L^4}^4+ \|\sqrt\rho \dot u\|_{L^2}^4\right).
\end{aligned}
\end{equation}

\underline{\bf Step 4. Proof of (\ref{3.3})}

\vskip 2mm

In view of (\ref{3.15}), (\ref{3.16}), (\ref{3.19})--(\ref{3.21}) and (\ref{3.31}), we derive from (\ref{3.13}) that
\begin{equation}
\begin{aligned}\label{3.32}
&\frac{\md}{\md t} \|\sqrt\rho \dot{u}\|_{L^2}^2 +
 \left(\mu\|\nabla^\bot\cdot \dot u\|_{L^2}^2+\nu\| \dot V\|_{L^2}^2 \right)\\
&\quad\leq -2\frac{\md}{\md t}\int  F   \partial_i u\cdot\nabla u_i\mdx  + C\left( 1+\|\nabla u\|_{L^4}^4+\|\sqrt\rho\dot u\|_{L^4}^4\right).
\end{aligned}
\end{equation}

It remains to deal with the first term on the right-hand side of (\ref{3.32}). Noting that
$$
 \partial_i u\cdot\nabla u_i=(\divg u)^2+2(\partial_1 u_2\partial_2 u_1-\partial_1 u_1\partial_2 u_2)=(\divg u)^2+2\nabla^\bot u_2\cdot\nabla u_1,
$$
we have
\begin{equation}\label{3.33}
2\int  F   \partial_i u\cdot\nabla u_i\mdx=2\int F(\divg u)^2\mdx+4\int F\nabla^\bot u_2\cdot\nabla u_1\mdx.
\end{equation}

Recalling that $\divg u=\nu^{-1}(F+P-\overline P)$ and noting that $\nu^{-1}\|F\|_{L^2}\leq C$ due to (\ref{3.1}) and (\ref{3.2}), by (\ref{2.1}) and (\ref{3.18}) we deduce
\begin{equation}\label{3.34}
\begin{aligned}
\left|2\int F(\divg u)^2\mdx\right| &=\frac{2}{\nu}\int F (F+P-\overline P)\divg u\mdx\\
&\leq C\left(\nu^{-1}\|F\|_{L^4}^2+\nu^{-1}\|F\|_{L^2}\right)\|\divg u\|_{L^2}\\
&\leq C\left(\nu^{-1}\|F\|_{L^2}\|\nabla F\|_{L^2}+\nu^{-1}\|F\|_{L^2}\right)\|\divg u\|_{L^2}\\
&\leq C+C\|\sqrt\rho \dot u\|_{L^2}.
\end{aligned}
\end{equation}

Next, since $\text{BMO}(\mathbb{T}^2)$ and $\mathcal{H}^1(\mathbb{T}^2)$ are dual spaces and  $\nabla^\bot\cdot (\nabla u_1)=\divg (\nabla^\bot u_2)=0$, by (\ref{2.7}) we infer from (\ref{3.2}) and (\ref{3.18}) that
\begin{equation}\label{3.35}
\begin{aligned}
\left|4\int F\nabla^\bot u_2\cdot\nabla u_1\mdx\right|&\leq C\|F\|_{\text{BMO}}\|\nabla^\bot u_2\cdot \nabla u_1\|_{\mathcal{H}^1}\\
&\leq C\|\nabla F\|_{L^2}\|\nabla u\|_{L^2}^2\leq C\|\sqrt\rho\dot u\|_{L^2}.
\end{aligned}
\end{equation}

Now, substituting (\ref{3.34}) and (\ref{3.35}) into (\ref{3.33}), we arrive at
\begin{equation}\label{3.36}
\left|2\int  F   \partial_i u\cdot\nabla u_i\mdx\right|\leq  C+C\|\sqrt\rho\dot u\|_{L^2}.
\end{equation}
Thus, multiplying (\ref{3.32}) by $t$, integrating it over $(0,T)$,  using (\ref{3.2}), (\ref{3.36}), Cauchy-Schwarz inequality and  Gronwall inequality, we obtain  the desired estimate stated in (\ref{3.3}).\hfill$\square$

\vskip 2mm

As an immediate consequence of Lemmas \ref{lem3.1} and \ref{lem3.2}, we have

\begin{lem}\label{lem3.3} Assume that the conditions of Proposition  \ref{pro1.2} hold. Then for any $2<q<\infty$ and $0<\varepsilon<1$,
\begin{equation}\label{3.37}
\begin{aligned}
&\int_0^T \left(\|\rho\dot{u}\|_{L^q}^{2-\varepsilon}+\|\nabla F\|_{L^q}^{2-\varepsilon}+\|\nabla(\nabla^\bot\cdot  u)\|_{L^q}^{2-\varepsilon}\right)\md t\\
&\quad+\int_0^T\left(\|F\|_{L^\infty}^{2-\varepsilon}+ \|\nabla^\bot\cdot u\|_{L^\infty}^{2-\varepsilon}+\nu^{2-\varepsilon}\|\divg u\|_{L^\infty}^{2-\varepsilon}\right)\md t\leq C(\varepsilon,T).
\end{aligned}
\end{equation}
\end{lem}
\pf First, it is easily derived from (\ref{3.17})  that for any $2\leq q<\infty$,
\begin{equation}\label{3.38}
\|\nabla F\|_{L^q}+\|\nabla(\nabla^\bot\cdot  u)\|_{L^q}\leq C \|\rho\dot u\|_{L^q}.
\end{equation}

For any $2<q<\infty$ and $0<\theta<1$, it follows from the interpolation inequality and (\ref{3.1}) that there exists a $\tilde q>q$, depending on $q$ and $\theta$, such that
$$
\|\rho\dot{u}\|_{L^q} \leq C \|\sqrt\rho\dot{u}\|_{L^2}^{\theta}\|\dot{u}\|_{L^{\tilde q}}^{1-\theta}\leq C \|\sqrt\rho\dot{u}\|_{L^2}^{\theta}\|\dot{u}\|_{L^{\tilde q}}^{1-\theta},
$$
which, combined with (\ref{2.4}), (\ref{2.5}) and (\ref{3.14}), shows (keeping in mind that $V=\divg u$)
\begin{align*}
\|\rho\dot{u}\|_{L^q} &\leq    C \|\sqrt\rho\dot{u}\|_{L^2}^{\theta}\left(\|\sqrt\rho\dot{u}\|_{L^2}+\|\nabla\dot u\|_{L^2}\right)^{1-\theta}\\
&\leq  C  \|\sqrt\rho\dot{u}\|_{L^2}^{\theta}\left(\|\sqrt\rho\dot{u}\|_{L^2}+\|\nabla^\bot\cdot\dot u\|_{L^2}+\|\divg \dot u\|_{L^2}\right)^{1-\theta}\\
&\leq  C  \|\sqrt\rho\dot{u}\|_{L^2}^{\theta}\left(\|\sqrt\rho\dot{u}\|_{L^2}+\|\nabla^\bot\cdot\dot u\|_{L^2}+\|\dot V\|_{L^2}+\|\nabla u\|_{L^4}^2\right)^{1-\theta} \\
&\leq C\left(\|\sqrt\rho\dot{u}\|_{L^2}+\|\nabla u\|_{L^4}^2\right)+C \|\sqrt\rho\dot{u}\|_{L^2}^{\theta}\left( \|\nabla^\bot\cdot\dot u\|_{L^2}+\|\dot V\|_{L^2}\right)^{1-\theta}.
\end{align*}
Thus, by virtue of (\ref{3.2}) and (\ref{3.3}) we deduce  that for any $0<\varepsilon<1$,
\begin{equation}\label{3.39}
\begin{aligned}
\int_0^T\|\rho\dot{u}\|_{L^q}^{2-\varepsilon}\md t&\leq C +C\int_0^T\left(\|\sqrt\rho\dot{u}\|_{L^2}^2+\|\nabla u\|_{L^4}^4\right)\md t\\
&\quad +C\int_0^T t^{-\frac{2-\varepsilon}{2}}\left(t \|\sqrt\rho\dot{u}\|_{L^2}^2\right)^{\frac{\theta(2-\varepsilon)}{2}}\left(t \|\nabla^\bot\cdot\dot u\|_{L^2}^2+t\|\dot V\|_{L^2}^2 \right)^{\frac{(1-\theta)(2-\varepsilon)}{2}}\md t\\
& \leq  C+C\sup_{0\leq t\leq T}\left(t \|\sqrt\rho\dot{u}\|_{L^2}^2\right)^{\frac{\theta(2-\varepsilon)}{2}}\left(\int_0^T t^{-\frac{2-\varepsilon}{\varepsilon+\theta(2-\varepsilon)}}\md t\right)^{\frac{\varepsilon+\theta(2-\varepsilon)}{2}}\\
&\qquad\qquad\qquad   \times\left(\int_0^T\left(t  \|\nabla^\bot\cdot\dot u\|_{L^2}^2+t\|\dot V\|_{L^2}^2 \right) \md t\right)^{\frac{(1-\theta)(2-\varepsilon)}{2}}\\
&\leq C,
\end{aligned}
\end{equation}
provided $\theta\in(0,1)$ is taken to be such that
$$
0<\frac{2(1-\varepsilon)}{2-\varepsilon}< \theta<1\quad \left(\text{and\ \ thus,}\quad   0<\frac{2-\varepsilon}{\varepsilon+\theta(2-\varepsilon)}<1\right).
$$
As a result, we infer from (\ref{2.2}), (\ref{3.2}), (\ref{3.18}) and (\ref{3.38}) that for any $2<q<\infty$ and $0<\varepsilon<1$,
\begin{equation}\label{3.40}
\begin{aligned}
&\int_0^T \left(\|F\|_{L^\infty}^{2-\varepsilon}+ \|\nabla^\bot\cdot u\|_{L^\infty}^{2-\varepsilon}\right)\md t\\
&\quad \leq C \int_0^T \left( \|F\|_{L^2}^{2-\varepsilon}+\|\nabla F\|_{L^q}^{2-\varepsilon}+\|\nabla^\bot\cdot u\|_{L^2}^{2-\varepsilon}+\|\nabla (\nabla^\bot\cdot u)\|_{L^q}^{2-\varepsilon}\right)\md t\\
&\quad \leq C+C\int_0^T \left(\|\sqrt\rho \dot u\|_{L^2}^{2-\varepsilon}+\|\rho\dot u\|_{L^q}^{2-\varepsilon} \right)\md t\leq C,
\end{aligned}
\end{equation}
and consequently,
\begin{equation*}
\nu^{2-\varepsilon}\int_0^T \|\divg u\|_{L^\infty}^{2-\varepsilon}\md t\leq C\int_0^T \left(\|F\|_{L^\infty}^{2-\varepsilon}+ \|P-\overline P\|_{L^\infty}^{2-\varepsilon}\right)\md t\leq C.
\end{equation*}
This, together (\ref{3.39}) and (\ref{3.40}), finishes the proof of (\ref{3.37}).\hfill$\square$

\vskip 2mm

\noindent{\bf  Proofs of Proposition \ref{pro1.2} and Theorem \ref{thm1.1}.} Based on Lemmas \ref{lem3.2} and \ref{lem3.3}, we immediately obtain Proposition \ref{pro1.2}. Then, combining Lemma \ref{lem3.1} and Proposition \ref{pro1.2}, we can prove the global existence theorem and the uniqueness result stated in Theorem \ref{thm1.1} in the same manner as that in \cite[Sections 4 and 5]{Danchin2023} .\hfill$\square$

\section{Convergence rate of the incompressible limit}\label{sec4}

With the help of the uniform-in-$\nu$ estimates established in Lemma \ref{lem3.1}--\ref{lem3.3}, the justification of the incompressible limit as $\nu\to\infty$ can be done in the same way as that in \cite{Danchin2023} and the details are omitted for simplicity. In this section, we focus on the convergence rate of the incompressible limit. To do this, we need the $L^q$-estimate  of the gradient of   density.

\begin{lem}\label{lem4.1}For any $q\in(2,\infty)$, assume that $\nabla\rho_0\in L^q$ holds. Then for any $0<\varepsilon<1$, there exists a positive constant $C$, depending on $\varepsilon$, $q$ and $T$, such that for any $\nu\geq\nu_0$,
\begin{equation}
\sup_{0\leq t\leq T}\|\nabla\rho(t)\|_{L^q}+\int_0^T\left(\nu^{2-\varepsilon}\|\nabla\divg u\|_{L^q}^{2-\varepsilon}+\|\nabla u\|_{W^{1,q}}^{2-\varepsilon}+\|\nabla u\|_{L^\infty} 
+\|\nabla u\|_{H^1}^2\right)\md t\leq C.\label{4.1}
\end{equation}
\end{lem}
\pf For any $2<q <\infty$, operating $\nabla$ to both sides of (\ref{1.1})$_1$ and multiplying it by $q|\nabla \rho|^q\nabla \rho$, by (\ref{1.10}) we obtain after integrating the resulting equation over $\mathbb{T}^2$ that for any $\nu\geq1$,
\begin{align*}
\frac{\md}{\md t}\|\nabla\rho\|_{L^q}^q&\leq C(q)\|\nabla u\|_{L^\infty}\|\nabla\rho\|_{L^q}^q+C(q)\|\nabla\divg u\|_{L^q}\|\nabla\rho\|_{L^q}^{q-1}\\
&\leq C(q)\|\nabla u\|_{L^\infty}\|\nabla\rho\|_{L^q}^q+C(q)\left(\|\nabla F\|_{L^q}+\|\nabla\rho\|_{L^q}\right)\|\nabla\rho\|_{L^q}^{q-1},
\end{align*}
which, combined with (\ref{3.38}), implies
\begin{equation}\label{4.2}
\frac{\md}{\md t}\|\nabla\rho\|_{L^q} \leq C(q)\left(1+\|\nabla u\|_{L^\infty}\right)\|\nabla\rho\|_{L^q}+C(q)\|\rho\dot u\|_{L^q}.
\end{equation}

By virtue of of (\ref{1.10}), (\ref{3.1}), (\ref{3.2})  and (\ref{3.38}), we see that
\begin{equation}\label{4.3}
\begin{aligned}
\|\nabla u\|_{W^{1,q}}&\leq  C\left(\|\nabla^\bot\cdot u\|_{W^{1,q}}+\|\divg u\|_{W^{1,q}}\right)\\
&\leq  C\left(\|\nabla^\bot\cdot u\|_{W^{1,q}}+\|F\|_{W^{1,q}}+\|\rho\|_{W^{1,q}}\right)\\
&\leq  C\left(1+\|\rho\dot u\|_{L^q}+ \|\nabla\rho\|_{L^q}\right).
\end{aligned}
\end{equation}

Using (\ref{3.2}) and (\ref{4.3}), we infer from (\ref{2.6}) that for any $2<p<\infty$,
\begin{equation}\label{4.4}
\begin{aligned}
\|\nabla u\|_{L^\infty}&\leq C\left(\|\divg u\|_{L^\infty}+\|\nabla^\bot\cdot u\|_{L^\infty}\right)\ln(e+\|\nabla^2 u\|_{L^q})+C\|\nabla u\|_{L^2}+C\\
&\leq C \left(1+\|\divg u\|_{L^\infty}+\|\nabla^\bot\cdot u\|_{L^\infty}\right)\ln \left(e+\|\rho\dot u\|_{L^q}+ \|\nabla\rho\|_{L^q}\right).
\end{aligned}
\end{equation}

Now, putting (\ref{4.4}) into (\ref{4.2}), we obtain
\begin{equation}\label{4.5}
\frac{\md}{\md t}\ln \left(e+\|\nabla\rho\|_{L^q}\right)\leq f(t) \ln \left(e+\|\nabla\rho\|_{L^q}\right)+g(t),
\end{equation}
where
$$
f(t):= C\left(1+\|\divg u\|_{L^\infty}+\|\nabla^\bot\cdot u\|_{L^\infty}\right)
$$
and
$$
g(t):= C \left(1+\|\rho\dot u\|_{L^q}+\|\divg u\|_{L^\infty}+\|\nabla^\bot\cdot u\|_{L^\infty}\right)\ln \left(e+\|\rho\dot u\|_{L^q}\right).
$$

It is easily deduced from (\ref{3.37}) that $f(t),g(t)\in L^1(0,T)$. Hence, an application of Gronwall's inequality to (\ref{4.5}) 
proves that  $\ln \left(e+\|\nabla\rho\|_{L^q}\right)$  is uniformly bounded in $\nu$ for any $t\in[0,T]$, which particularly yields an uniform-in-$\nu$ bound of  $\|\nabla\rho\|_{L^q}$ on  $[0,T]$. As an easy result, the estimates for $\|\nabla u\|_{L^\infty}$, $\nu^{2-\varepsilon}\|\nabla\divg u\|_{L^q}^{2-\varepsilon}$  and $\|\nabla u\|_{H^1}^2$ follow   from (\ref{3.2}),  (\ref{3.37})  and  (\ref{4.4}) immediately. \hfill$\square$

\vskip 2mm

Based upon the global regularity of $(\eta,v)$ given in  (\ref{2.8}) and the uniform estimates established in Lemmas \ref{lem3.1}--\ref{lem3.3} and \ref{lem4.1}, we are now ready to derive the convergence rate of the incompressible limit as $\nu\to\infty$.

\vskip 2mm

\noindent{\bf Proof of Theorem \ref{thm1.2}.} Assume that $(\rho,u)$ and $(\eta,v)$ are the solutions of (\ref{1.1}) and (\ref{1.14}) with the same initial data $(\rho_0,u_0)$, respectively. To begin, we recall that $\mathcal{P} u:=(\mathbb{Id}+\nabla(-\Delta)^{-1}\divg )u$, $\mathcal{Q}u:=-\nabla(-\Delta)^{-1}\divg u$ and
\begin{equation}\label{4.6}
u= \mathcal{P} u+\mathcal{Q}u\quad{\text{with}}\quad  \divg (\mathcal{P}u)=0,\quad \nabla^\bot\cdot(\mathcal{Q}u)=0.
\end{equation}

\underline{\bf Step I. The $t$-growth and singular $t$-weighted estimates} 

\vskip 2mm

Let  $\tilde\rho=\tilde\rho(x,t)$  be the ``incompressible" density related to the divergence-free part of $u$ and determined by
\begin{equation}\label{4.7}
\tilde\rho_t+\mathcal{P}u\cdot\nabla\tilde\rho=0,\quad \tilde\rho|_{t=0}=\rho_0.
\end{equation}
Thanks to (\ref{4.1}), it is easy to get that for any $0\leq t\leq T$,
\begin{equation}
\begin{aligned}
\|\nabla\tilde\rho(t)\|_{L^q}&\leq C\|\nabla\rho_0\|_{L^q}\exp\left\{\int_0^t\|\nabla (\mathcal{P}u)\|_{L^\infty}\md\tau\right\}\\
&\leq C\|\nabla\rho_0\|_{L^q}\exp\left\{\int_0^t\|\nabla  u\|_{W^{1,q}}\md\tau\right\}\leq C.
\end{aligned}\label{4.8}
\end{equation}

The difference of $\rho-\eta$ will be decomposed into two parts,
$$\rho-\eta=(\rho-\tilde\rho)+(\tilde\rho-\eta):=\varphi+\phi\quad\text{with}\quad \varphi:=\rho-\tilde\rho,\quad \phi:=\tilde\rho-\eta.$$ 
Then, it follows from (\ref{1.1})$_1$,   (\ref{4.6}) and (\ref{4.7}) that
\begin{equation}\label{4.9}
\varphi_t+ \mathcal{P}u\cdot\nabla \varphi=-\divg(\rho\mathcal{Q}u),\quad \varphi|_{t=0}=0
\end{equation}
and
\begin{equation}\label{4.10}
\phi_t+v\cdot\nabla \phi=-(\mathcal{P}u-v)\cdot\nabla\tilde\rho,\quad \phi|_{t=0}=0.
\end{equation}

Next, we aim to derive some $t$-growth and singular $t$-weighted estimates for $(\varphi,\phi)$, which will be used to overcome the difficulties induced by the presence of vacuum and to deal with the discrepancy of the velocities. We begin with the $t$-growth estimates of  $(\varphi,\phi)$.  Multiplying (\ref{4.9})  by $\varphi$   in $L^2$ and using the fact that $\divg (\mathcal{P}u)=\nabla^\bot\cdot(\mathcal{Q}u)=0$, by (\ref{2.1}), (\ref{2.5}), (\ref{3.2}) and (\ref{4.1}) we deduce that   for $\frac{1}{q}+\frac{1}{q^*}=\frac{1}{2}$ with $2<q<\infty$,
\begin{equation}\label{4.11}
\begin{aligned}
 \|\varphi(t)\|_{L^2}&\leq C\int_0^t\left( \|\rho\divg (\mathcal{Q}u)\|_{L^2}+\|\mathcal{Q}u\cdot\nabla \rho\|_{L^2}\right)\md\tau \\
&\leq C\int_0^t\left(\|\divg u\|_{L^2}+\|\mathcal{Q}u\|_{L^{q^*}}\|\nabla\rho\|_{L^q} \right)\md \tau\\
&\leq C\int_0^t\left(\|\divg u\|_{L^2}+\|\nabla \mathcal{Q}u\|_{L^2} \right)\md \tau\\
&\leq  C\int_0^t \|\divg u\|_{L^2} \md \tau\leq C\nu^{-\frac{1}{2}} t,
\end{aligned}
\end{equation}
where we have also used the identity  $\divg (\mathcal{Q}u)=\divg u$ and the Poincar${\rm\acute{e}}$ inequality due to the fact that the mean of $\mathcal{Q}u$ is zero. 

Analogously, we have by  multiplying (\ref{4.10})  by $\phi$   in $L^2$ that
\begin{equation}\label{4.12}
\frac{\md}{\md t}\|\phi(t)\|_{L^2}^2 \leq \||\mathcal{P}u-v||\nabla\tilde\rho|\|_{L^2}\|\phi\|_{L^2}\leq C \|\mathcal{P}u-v\|_{L^{q^*}}\|\nabla\tilde\rho\|_{L^q}\|\phi\|_{L^2},
\end{equation}
which, combined with (\ref{2.4}), (\ref{2.8}), (\ref{3.2})  and (\ref{4.8}), gives
\begin{equation}\label{4.13}
\begin{aligned}
 \|\phi(t)\|_{L^2} 
&\leq C\int_0^t\left( \|u\|_{L^{q^*}}+\| v\|_{L^{q^*}} \right)\md \tau\\
&\leq C\int_0^t\left(\|\sqrt \rho u\|_{L^{2}}+\|\nabla u\|_{L^2}+\|\sqrt\eta v\|_{L^2}+\|\nabla v\|_{L^2} \right)\md \tau\\
&\leq Ct.
\end{aligned}
\end{equation}

From (\ref{2.4}), (\ref{4.8}) and (\ref{4.12}), we also derive that
\begin{equation}\label{4.14}
\begin{aligned}
\frac{\md}{\md t}\|\phi(t)\|_{L^2}^2 &\leq  C \|\mathcal{P}u-v\|_{L^{q^*}}\|\nabla\tilde\rho\|_{L^q} \|\phi\|_{L^2}\\
&\leq  C \left(\|\sqrt\eta(\mathcal{P}u-v)\|_{L^2} +\|\nabla(\mathcal{P}u-v) \|_{L^2}\right)\|\phi\|_{L^2}\\
&\leq  C \left(\|\sqrt\eta(\mathcal{P}u-v)\|_{L^2}^2+\|\phi\|_{L^2}^2\right)+C_1\|\nabla(\mathcal{P}u-v) \|_{L^2}^2.
\end{aligned}
\end{equation}

We proceed to prove the singular $t$-weighted estimate of $\phi$.  Multiplying (\ref{4.10})  by $\frac{\phi}{t}$  in $L^2$ and integrating by parts, by (\ref{2.4}) and (\ref{4.8}) we obtain 
\begin{align*}
\frac{1}{2}\frac{\md}{\md t}\left(\frac{\|\phi\|_{L^2}^2}{t} \right)+\frac{1}{2}  \frac{\|\phi\|_{L^2}^2}{t^2} 
&\leq  t^{-1} \|\mathcal{P}u-v\|_{L^{q^*}}\|\nabla\tilde\rho\|_{L^q}\|\phi\|_{L^2} \\
&\leq \frac{ \|\phi\|_{L^2}^2}{4t^2} +C\left(\|\sqrt\eta(\mathcal{P}u-v) \|_{L^2}^2+\|\nabla (\mathcal{P}u-v)\|_{L^2}^2\right),
\end{align*}
and hence,
\begin{equation}
 \frac{\md}{\md t}\left(\frac{\|\phi\|_{L^2}^2}{t} \right)+  \frac{\|\phi\|_{L^2}^2}{2t^2} 
\leq C\|\sqrt\eta(\mathcal{P}u-v) \|_{L^2}^2+C_2\|\nabla (\mathcal{P}u-v)\|_{L^2}^2 .\label{4.15}
\end{equation}

\vskip 2mm

\underline{\bf Step II. The difference of velocities}

\vskip 2mm

First, it follows from (\ref{1.1})$_2$ and (\ref{4.6}) that the divergence-free part of the velocity $ \mathcal{P}u$ satisfies
\begin{equation}\label{4.16}
\begin{aligned}
&\tilde\rho (\mathcal{P}u)_t+\tilde\rho \mathcal{P}u\cdot\nabla(\mathcal{P}u)-\mu \Delta (\mathcal{P}u)+\nabla \Lambda \\
&\quad=  -\varphi(\mathcal{P}u)_t -\rho (\mathcal{Q}u)_t -\varphi \mathcal{P}u\cdot\nabla (\mathcal{P}u)-\rho\mathcal{Q}u\cdot\nabla(\mathcal{P}u)-\rho u\cdot\nabla(\mathcal{Q}u),
\end{aligned}
\end{equation}
where $\Lambda:=-\nu\divg u+( P(\rho)-\overline P)$. Then, subtracting (\ref{1.14})$_2$ from (\ref{4.16}) gives
\begin{equation}\label{4.17}
\begin{aligned}
&\eta (\mathcal{P}u-v)_t+\eta v\cdot\nabla(\mathcal{P}u-v)-\mu\Delta(\mathcal{P}u-v)+\nabla (\Lambda-\Pi)\\
&\quad= -\phi(\mathcal{P}u)_t-\varphi(\mathcal{P}u)_t-\phi\mathcal{P}u\cdot\nabla(\mathcal{P}u)-\varphi \mathcal{P}u\cdot\nabla (\mathcal{P}u)\\
&\qquad -\eta (\mathcal{P}u-v)\cdot\nabla(\mathcal{P}u)-\rho\mathcal{Q}u\cdot\nabla(\mathcal{P}u)-\rho u\cdot\nabla(\mathcal{Q}u)-\rho (\mathcal{Q}u)_t,
\end{aligned}
\end{equation}
which is supplemented with vanishing initial condition $(\mathcal{P}u-v)|_{t=0}=0$. 

Next, in terms of (\ref{1.14})$_1$ and the divergence-free condition $\divg (\mathcal{P}u-v)=0$, we obtain after multiplying (\ref{4.17}) by $(\mathcal{P}u-v)$ in $L^2$ and integrating by parts that
\begin{equation}\label{4.18}
\begin{aligned}
&\frac{1}{2}\frac{\md}{\md t}\int \eta |\mathcal{P}u-v|^2\mdx+\mu\int |\nabla(\mathcal{P}u-v)|^2\mdx \\
&\quad= -\int \phi(\mathcal{P}u)_t\cdot  (\mathcal{P}u-v)\mdx-\int \varphi(\mathcal{P}u)_t \cdot (\mathcal{P}u-v)\mdx\\
&\qquad -\int \phi\mathcal{P}u\cdot\nabla(\mathcal{P}u) \cdot (\mathcal{P}u-v)\mdx -\int \varphi \mathcal{P}u\cdot\nabla (\mathcal{P}u)\cdot  (\mathcal{P}u-v)\mdx\\
&\qquad-\int \eta (\mathcal{P}u-v)\cdot\nabla(\mathcal{P}u)  \cdot (\mathcal{P}u-v)\mdx-\int \rho\mathcal{Q}u\cdot\nabla(\mathcal{P}u) \cdot (\mathcal{P}u-v)\mdx  \\
&\qquad  -\int \rho u\cdot\nabla(\mathcal{Q}u) \cdot (\mathcal{P}u-v)\mdx-\int \rho (\mathcal{Q}u)_t  \cdot (\mathcal{P}u-v)\mdx \\
&\quad :=M_1+\ldots+M_8.
\end{aligned}
\end{equation}

We are now in a position of considering each term on the right-hand side of (\ref{4.18}). First, by virtue of  (\ref{2.4}) we deduce that for any $2<r<\infty$ and $\frac{1}{r}+\frac{1}{r^*}=\frac{1}{2}$,
\begin{equation}\label{4.19}
\begin{aligned}
M_1&\leq C \frac{\|\phi\|_{L^2}}{\sqrt t} \left(\sqrt t\left\| (\mathcal{P}u)_t\right\|_{L^{r^*}}\right)\left\| \mathcal{P}u-v \right\|_{L^{r}} \\
&\leq C\frac{\|\phi\|_{L^2}}{\sqrt t} \left(\sqrt t\left\| \mathcal{P}u_t\right\|_{L^{r^*}}\right)\left(\left\|\sqrt\eta(\mathcal{P}u-v)\right\|_{L^2}+\left\|\nabla(\mathcal{P}u-v)\right\|_{L^2}\right) \\
&\leq \frac{\mu}{20}\left \|\nabla(\mathcal{P}u-v)\right\|_{L^2}^2+C \left\|\sqrt\eta(\mathcal{P}u-v)\right\|_{L^2}^2+ C\frac{\|\phi\|_{L^2}^2}{  t} \left( t\left\| \mathcal{P}u_t\right\|_{L^{r^*}}^2\right). 
\end{aligned}
\end{equation}

In a similar manner, by (\ref{4.11}) we have
\begin{equation}\label{4.20}
\begin{aligned}
M_2&\leq  C \frac{\|\varphi\|_{L^2}}{\sqrt t} \left(\sqrt t\left\| (\mathcal{P}u)_t\right\|_{L^{r^*}}\right)\left\| \mathcal{P}u-v \right\|_{L^{r}} \\
&\leq C\frac{\|\varphi\|_{L^2}}{\sqrt t} \left(\sqrt t\left\| \mathcal{P}u_t\right\|_{L^{r^*}}\right)\left(\left\|\sqrt\eta(\mathcal{P}u-v)\right\|_{L^2}+\left\|\nabla(\mathcal{P}u-v)\right\|_{L^2}\right) \\
&\leq \frac{\mu}{20}\left \|\nabla(\mathcal{P}u-v)\right\|_{L^2}^2+C \left\|\sqrt\eta(\mathcal{P}u-v)\right\|_{L^2}^2+ C\nu^{-1} \left( t\left\| \mathcal{P}u_t\right\|_{L^{r^*}}^2\right). 
\end{aligned}
\end{equation}

Due to (\ref{2.4}),   (\ref{3.1}) and (\ref{3.2}), it holds that
\begin{equation}\label{4.21}
\|\mathcal{P}u\|_{L^\alpha}\leq C\|u\|_{L^\alpha}\leq  C(\alpha)\left(\|\sqrt\rho u\|_{L^2}+\|\nabla u\|_{L^2}\right)\leq C(\alpha),\quad\forall \ 2\leq \alpha<\infty,
\end{equation}
and hence,
\begin{equation}\label{4.22}
\begin{aligned}
M_3&\leq  \| \phi  \|_{L^2}\|\mathcal{P}u\|_{L^8}\left\|\nabla (\mathcal{P}u)\right\|_{L^4}\left\| \mathcal{P}u-v \right\|_{L^{8}} \\
&\leq \left\| \phi \right\|_{L^2} \|\nabla u\|_{L^4} \left(\left\|\sqrt\eta(\mathcal{P}u-v)\right\|_{L^2}+\left\|\nabla(\mathcal{P}u-v)\right\|_{L^2}\right) \\
&\leq \frac{\mu}{20}\left \|\nabla(\mathcal{P}u-v)\right\|_{L^2}^2+C \left\|\sqrt\eta(\mathcal{P}u-v)\right\|_{L^2}^2+ C \|\phi\|_{L^2}^2 \|\nabla u\|_{L^4}^2. 
\end{aligned}
\end{equation}
Analogously, by (\ref{4.11}) and (\ref{4.21}) we have
\begin{equation}\label{4.23}
\begin{aligned}
M_4&\leq  \| \varphi  \|_{L^2}\|\mathcal{P}u\|_{L^8}\left\|\nabla (\mathcal{P}u)\right\|_{L^4}\left\| \mathcal{P}u-v \right\|_{L^{8}} \\
&\leq \left\| \varphi \right\|_{L^2} \|\nabla u\|_{L^4} \left(\left\|\sqrt\eta(\mathcal{P}u-v)\right\|_{L^2}+\left\|\nabla(\mathcal{P}u-v)\right\|_{L^2}\right)\\
&\leq \frac{\mu}{20}\left \|\nabla(\mathcal{P}u-v)\right\|_{L^2}^2+C \left\|\sqrt\eta(\mathcal{P}u-v)\right\|_{L^2}^2+ C \nu^{-1} \|\nabla u\|_{L^4}^2. 
\end{aligned}
\end{equation}

It is easily seen from (\ref{2.1}), (\ref{2.4}), (\ref{2.5}), (\ref{2.8}), (\ref{3.1}) and (\ref{3.2}) that
\begin{equation}
M_5\leq C\left\|\nabla(\mathcal{P} u)\right\|_{L^\infty} \left\|\sqrt\eta(\mathcal{P}u-v)\right\|_{L^2}^2\leq C\left\|\nabla u\right\|_{W^{1,q}} \left\|\sqrt\eta(\mathcal{P}u-v)\right\|_{L^2}^2.\label{4.24}
\end{equation}
and
\begin{equation}\label{4.25}
\begin{aligned}
M_6+M_7&\leq C\left(\|\mathcal{Q}u\|_{L^4}\|\nabla \mathcal{P}u\|_{L^2} +\| u\|_{L^4}\|\nabla \mathcal{Q}u\|_{L^2}\right)\|\mathcal{P}u-v\|_{L^4}\\
&\leq C\|u\|_{H^1}\|\nabla(\mathcal{Q}u)\|_{L^2}\left(\|\eta(\mathcal{P}u-v)\|_{L^2}+\|\nabla(\mathcal{P}u-v)\|_{L^2}\right)\\
&\leq C \|\divg u \|_{L^2}\left(\|\eta(\mathcal{P}u-v)\|_{L^2}+\|\nabla(\mathcal{P}u-v)\|_{L^2}\right)\\
&\leq   \frac{\mu}{20}\left \|\nabla(\mathcal{P}u-v)\right\|_{L^2}^2+C \left\|\sqrt\eta(\mathcal{P}u-v)\right\|_{L^2}^2+C\nu^{-1}.
\end{aligned}
\end{equation}

The handling of the last term $M_8$ needs more works. Recalling that $\rho=\varphi+\phi+\eta$, we have
\begin{equation}\label{4.26}
M_8= -\int \left(\phi+\varphi\right) (\mathcal{Q}u)_t  \cdot (\mathcal{P}u-v)\mdx  -\int \eta (\mathcal{Q}u)_t \cdot (\mathcal{P}u-v)\mdx:=M_{8,1}+M_{8,2}.
\end{equation}

The  first term on the right-hand side of (\ref{4.26}) can be bounded in a manner similar to that in (\ref{4.19}) and (\ref{4.20}) by
\begin{equation}\label{4.27}
\begin{aligned}
M_{8,1}&\leq  C \frac{\|(\phi,\varphi)\|_{L^2}}{\sqrt t} \left(\sqrt t\left\| (\mathcal{Q}u)_t\right\|_{L^{r^*}}\right)\left\| \mathcal{P}u-v \right\|_{L^{r}}\\
&\leq C\frac{\|(\phi,\varphi)\|_{L^2}}{\sqrt t} \left(\sqrt t\left\| \mathcal{Q}u_t\right\|_{L^{r^*}}\right)\left(\left\|\sqrt\eta(\mathcal{P}u-v)\right\|_{L^2}+\left\|\nabla(\mathcal{P}u-v)\right\|_{L^2}\right)\\
&\leq \frac{\mu}{20}\left \|\nabla(\mathcal{P}u-v)\right\|_{L^2}^2+C \left\|\sqrt\eta(\mathcal{P}u-v)\right\|_{L^2}^2+ C \left( \frac{\|\phi\|_{L^2}^2}{  t} + \nu^{-1}\right) \left( t\left\| \mathcal{Q}u_t\right\|_{L^{r^*}}^2\right),
\end{aligned}
\end{equation}
where $2<r<\infty$ and $\frac{1}{r}+\frac{1}{r^*}=\frac{1}{2}$. To deal with the second term $M_{8,2}$, by (\ref{1.14})$_1$ we first rewrite it in the form:
\begin{equation}\label{4.28}
\begin{aligned}
M_{8,2}&=-\frac{\md}{\md t}\int \eta (\mathcal{Q}u) \cdot (\mathcal{P}u-v)\mdx+\int(\mathcal{Q}u)\cdot \left[ \eta_t  (\mathcal{P}u-v) +\eta(\mathcal{P}u_t-v_t)\right]\mdx\\
&=-\frac{\md}{\md t}\int \eta (\mathcal{Q}u) \cdot (\mathcal{P}u-v)\mdx-\int v\cdot\nabla\eta (\mathcal{Q}u) \cdot (\mathcal{P}u-v)\mdx\\
&\quad +\int \eta  (\mathcal{Q}u)  \cdot (\mathcal{P}u_t-v_t)\mdx\\
&:= -\frac{\md}{\md t}\int \eta (\mathcal{Q}u)  (\mathcal{P}u-v)\mdx+M_{8,2}^1+M_{8,2}^2.
\end{aligned}
\end{equation}

As before, let $\frac{1}{q}+\frac{1}{q^*}=\frac{1}{2}$ and $\frac{1}{r}+\frac{1}{r^*}=\frac{1}{2}$ for $2<q,r<\infty$. Then, using (\ref{2.4}), (\ref{2.5}), (\ref{2.8}), (\ref{3.2}) and the Poincar${\rm {\acute e}}$ inequality that (because of $\overline{\mathcal{Q}u}=0$), we deduce
\begin{equation}\label{4.29}
\begin{aligned}
M_{8,2}^1&\leq C\|v\|_{L^{q^*}}\|\nabla\eta\|_{L^q}\|\mathcal{Q}u\|_{L^{r^*}}\|\mathcal{P}u-v\|_{L^r}\\
&\leq C\left(\|\sqrt\eta v\|_{L^2}+\|\nabla v\|_{L^2}\right)\|\nabla\mathcal{Q}u\|_{L^2}\left(\|\sqrt\eta (\mathcal{P}u-v)\|_{L^2}+\|\nabla(\mathcal{P}u-v)\|\right)\\
&\leq C \|\divg u\|_{L^2}\left(\|\sqrt\eta (\mathcal{P}u-v)\|_{L^2}+\|\nabla(\mathcal{P}u-v)\|\right)\\
&\leq  \frac{\mu}{20}\left \|\nabla(\mathcal{P}u-v)\right\|_{L^2}^2+C \left\|\sqrt\eta(\mathcal{P}u-v)\right\|_{L^2}^2+C\nu^{-1}.
\end{aligned}
\end{equation}

Noting that $u_t=\mathcal{P}u_t+\mathcal{Q}u_t$ and $\eta=\eta^{\frac{1}{2}} (\rho-\phi-\varphi)^{\frac{1}{2}}$, by (\ref{1.14})$_1$ we have
\begin{equation}\label{4.30}
\begin{aligned}
M_{8,2}^2&=-\int \eta  (\mathcal{Q}u)  \mathcal{Q}u_t\mdx+\int \eta  (\mathcal{Q}u)  (u_t-v_t)\mdx\\
&\leq -\frac{1}{2}\frac{\md}{\md t}\int \eta |\mathcal{Q}u|^2\mdx+C\int |v||\nabla\eta| |\mathcal{Q}u|^2\mdx \\
&\quad + C\int \eta^{\frac{1}{2}}
\rho^{\frac{1}{2}} |\mathcal{Q}u| | u_t-v_t|\mdx
+ C\int \eta^{\frac{1}{2}}
(|\phi|+|\varphi|)^{\frac{1}{2}} |\mathcal{Q}u| | u_t-v_t|\mdx\\
&:= -\frac{1}{2}\frac{\md}{\md t}\int \eta |\mathcal{Q}u|^2\mdx+M_{8,2}^{2,1}+M_{8,2}^{2,2}+M_{8,2}^{2,3}.
\end{aligned}
\end{equation}

In terms of (\ref{2.4}), (\ref{2.5}), (\ref{2.8}), (\ref{3.2}) and the Poincar${\rm {\acute e}}$ inequality, we  find
\begin{equation}\label{4.31}
\begin{aligned}
M_{8,2}^{2,1}&\leq C\|\nabla\eta\|_{L^q}\|v\|_{L^q}\|\mathcal{Q}u\|_{L^{q^*}}^2\\
& \leq C\|\nabla(\mathcal{Q}u)\|_{L^2}^2\leq C\|\divg u\|_{L^2}^2\leq C\nu^{-1},
\end{aligned}
\end{equation}
and 
\begin{equation}\label{4.32}
\begin{aligned}
M_{8,2}^{2,2}&\leq C\|\mathcal{Q}u\|_{L^2}\left(\|\sqrt\rho u_t\|_{L^2}+\|\sqrt\eta v_t\|_{L^2}\right)\\
&\leq C\|\divg u\|_{L^2}\left(\|\sqrt\rho u_t\|_{L^2}+\|\sqrt\eta v_t\|_{L^2}\right)\\
&\leq C\nu^{-\frac{1}{2}}\left(\|\sqrt\rho u_t\|_{L^2}+\|\sqrt\eta v_t\|_{L^2}\right).
\end{aligned}
\end{equation}
Analogously, by  (\ref{4.11}) we obtain
\begin{equation}\label{4.33}
\begin{aligned}
M_{8,2}^{2,3}&\leq  \left(\|\phi\|_{L^2}^{\frac{1}{2}}+\|\varphi\|_{L^2}^{\frac{1}{2}}\right)\|\mathcal{Q}u\|_{L^2}\left(\|  u_t\|_{L^4}+\|  v_t\|_{L^4}\right)\\
&\leq C\|\divg u\|_{L^2}\left(\|\phi\|_{L^2}^{\frac{1}{2}}+\|\varphi\|_{L^2}^{\frac{1}{2}}\right)  \left(\|  u_t\|_{L^4}+\|  v_t\|_{L^4}\right)\\
&\leq C\nu^{-\frac{1}{2}}\left(t^{-\frac{1}{2}}\|\phi\|_{L^2}^{\frac{1}{2}}+\nu^{-\frac{1}{4}}\right)  \left(\sqrt t\|  u_t\|_{L^4}+\sqrt t\|  v_t\|_{L^4}\right).
\end{aligned}
\end{equation}

Inserting (\ref{4.31})--(\ref{4.33}) into (\ref{4.30}) gives
\begin{equation}\label{4.34}
\begin{aligned}
M_{8,2}^2&\leq  -\frac{1}{2}\frac{\md}{\md t}\int \eta |\mathcal{Q}u|^2\mdx + C\nu^{-\frac{1}{2}}\left(\|\sqrt\rho u_t\|_{L^2}+\|\sqrt\eta v_t\|_{L^2}\right)\\
&\quad + C\nu^{-\frac{1}{2}}\left(t^{-\frac{1}{2}}\|\phi\|_{L^2}^{\frac{1}{2}}+\nu^{-\frac{1}{4}}\right)  \left(\sqrt t\|  u_t\|_{L^4}+\sqrt t\|  v_t\|_{L^4}\right)+C\nu^{-1}.
\end{aligned}
\end{equation}
Hence, plugging (\ref{4.29}) and (\ref{4.34}) into (\ref{4.28}), we deduce
\begin{equation}\label{4.35}
\begin{aligned}
M_{8,2}&\leq   \frac{\mu}{20}\left \|\nabla(\mathcal{P}u-v)\right\|_{L^2}^2-\frac{\md}{\md t}\int \left(\eta (\mathcal{Q}u)  (\mathcal{P}u-v) +\frac{1}{2}  \eta |\mathcal{Q}u|^2 \right)\mdx\\
&\quad +C \left\|\sqrt\eta(\mathcal{P}u-v)\right\|_{L^2}^2+ C\nu^{-\frac{1}{2}}\left(\|\sqrt\rho u_t\|_{L^2}+\|\sqrt\eta v_t\|_{L^2}\right)\\
&\quad + C\nu^{-\frac{1}{2}}\left(t^{-\frac{1}{2}}\|\phi\|_{L^2}^{\frac{1}{2}}+\nu^{-\frac{1}{4}}\right)  \left(\sqrt t\|  u_t\|_{L^4}+\sqrt t\|  v_t\|_{L^4}\right)+C\nu^{-1}.
\end{aligned}
\end{equation}

As a result of  (\ref{4.27}) and (\ref{4.35}), we have from  (\ref{4.26}) that
\begin{equation}\label{4.36}
\begin{aligned}
M_{8}&\leq   \frac{\mu}{10}\left \|\nabla(\mathcal{P}u-v)\right\|_{L^2}^2-\frac{\md}{\md t}\int \left(\eta (\mathcal{Q}u)  (\mathcal{P}u-v) +\frac{1}{2}  \eta |\mathcal{Q}u|^2 \right)\mdx\\
&\quad +C \left\|\sqrt\eta(\mathcal{P}u-v)\right\|_{L^2}^2+ C\nu^{-\frac{1}{2}}\left(\|\sqrt\rho u_t\|_{L^2}+\|\sqrt\eta v_t\|_{L^2}\right)\\
&\quad + C\nu^{-\frac{1}{2}}\left(t^{-\frac{1}{2}}\|\phi\|_{L^2}^{\frac{1}{2}}+\nu^{-\frac{1}{4}}\right)  \left(\sqrt t\|  u_t\|_{L^4}+\sqrt t\|  v_t\|_{L^4}\right)\\
&\quad + C \left( t^{-1} \|\phi\|_{L^2}^2  + \nu^{-1}\right) \left( t\left\| \mathcal{Q}u_t\right\|_{L^{r^*}}^2\right)+C\nu^{-1}.
\end{aligned}
\end{equation}

Based upon (\ref{4.19}), (\ref{4.20}), (\ref{4.22})--(\ref{4.25}) and (\ref{4.36}), we infer from (\ref{4.18}) that
\begin{equation}\label{4.37}
\begin{aligned}
& \frac{\md}{\md t}\int \eta |\mathcal{P}u-v|^2\mdx+\mu\int |\nabla(\mathcal{P}u-v)|^2\mdx \\
&\quad \leq -\frac{\md}{\md t}\int\left(2 \eta (\mathcal{Q}u)  (\mathcal{P}u-v)+ \eta |\mathcal{Q}u|^2\right)\mdx +C\nu^{-1}  \\
&\qquad+C \left(1+\|\nabla u\|_{W^{1,q}}\right) \left\|\sqrt\eta(\mathcal{P}u-v)\right\|_{L^2}^2+C\|\phi\|_{L^2}^2\|\nabla u\|_{L^4}^2\\
&\qquad + C  t^{-1}\|\phi\|_{L^2}^2  \left( t\left\| \mathcal{P}u_t\right\|_{L^{r^*}}^2+ t\left\| \mathcal{Q}u_t\right\|_{L^{r^*}}^2\right)\\
&\qquad + C\nu^{-\frac{1}{2}}\left(t^{-\frac{1}{2}} \|\phi\|_{L^2}^{\frac{1}{2}}\right)  \left(\sqrt t\|  u_t\|_{L^4}+\sqrt t\|  v_t\|_{L^4}\right)\\
&\qquad  + C  \nu^{-1} \left(\|\nabla u\|_{L^4}^2+ t\left\| \mathcal{P}u_t\right\|_{L^{r^*}}^2+ t\left\| \mathcal{Q}u_t\right\|_{L^{r^*}}^2\right)\\
&\qquad+ C\nu^{-\frac{1}{2}}\left(\|\sqrt\rho u_t\|_{L^2}+\|\sqrt\eta v_t\|_{L^2}\right)+ C\nu^{-\frac{3}{4}}  \left(\sqrt t\|  u_t\|_{L^4}+\sqrt t\|  v_t\|_{L^4}\right).
\end{aligned}
\end{equation}

Let $\beta:=\frac{\mu}{4(C_1+C_2)}$ with $C_1$, $C_2$ being the positive numbers in (\ref{4.14}) and (\ref{4.15}). It is obvious that the sixth term on the right-hand side of (\ref{4.37}) can be bounded by
\begin{align*}
&C\nu^{-\frac{1}{2}}\left(t^{-\frac{1}{2}} \|\phi\|_{L^2}^{\frac{1}{2}}\right)  \left(\sqrt t\|  u_t\|_{L^4}+\sqrt t\|  v_t\|_{L^4}\right)\\
&\quad \leq \frac{\beta}{4}\frac{\|\phi\|_{L^2}^2}{t^2}+C(\beta)\nu^{-\frac{2}{3}} \left(\sqrt t\|  u_t\|_{L^4}+\sqrt t\|  v_t\|_{L^4}\right)^{\frac{4}{3}}.
\end{align*}
Thus, multiplying both (\ref{4.14}) and (\ref{4.15}) by   $\beta=\frac{\mu}{4(C_1+C_2)}$ and  adding them to (\ref{4.37}), we  obtain
\begin{equation}\label{4.38}
\begin{aligned}
& \frac{\md}{\md t} \left(\|\sqrt\eta (\mathcal{P}u-v)\|_{L^2}^2+\beta\left(\|\phi\|_{L^2}^2+\frac{\|\phi\|_{L^2}^2}{t}\right)\right)\\
&\qquad+\frac{1}{2}\left(\mu\|\nabla(\mathcal{P}u-v)\|_{L^2}^2+\frac{\beta}{2} \frac{\|\phi\|_{L^2}^2}{t^2}\right)\\
&\quad \leq -\frac{\md}{\md t}\int\left(2 \eta (\mathcal{Q}u)  (\mathcal{P}u-v)+ \eta |\mathcal{Q}u|^2\right)\mdx +C\nu^{-1}   \\
&\qquad +\left(\|\sqrt\eta (\mathcal{P}u-v)\|_{L^2}^2+ \|\phi\|_{L^2}^2+\frac{\|\phi\|_{L^2}^2}{t}\right) H(t) +K(t),
\end{aligned}
\end{equation}
where
$$
H(t):=C \left(1+\|\nabla u\|_{W^{1,q}}+\|\nabla u\|_{L^4}^2+t\left\| \mathcal{P}u_t\right\|_{L^{r^*}}^2+ t\left\| \mathcal{Q}u_t\right\|_{L^{r^*}}^2\right) 
$$
and
\begin{align*}
K(t)&:=C  \nu^{-1} \left(\|\nabla u\|_{L^4}^2+ t\left\| \mathcal{P}u_t\right\|_{L^{r^*}}^2+ t\left\| \mathcal{Q}u_t\right\|_{L^{r^*}}^2\right) + C\nu^{-\frac{1}{2}}\left(\|\sqrt\rho u_t\|_{L^2}+\|\sqrt\eta v_t\|_{L^2}\right)\\
&\quad+ C\nu^{-\frac{3}{4}}  \left(\sqrt t\|  u_t\|_{L^4}+\sqrt t\|  v_t\|_{L^4}\right) +C \nu^{-\frac{2}{3}} \left(\sqrt t\|  u_t\|_{L^4}+\sqrt t\|  v_t\|_{L^4}\right)^{\frac{4}{3}}.
\end{align*}
Here, the positive constant $C$ may depend on $\beta$.

\vskip 2mm

\underline{\bf Step III. The convergence rate of  incompressible limit}

\vskip 2mm

To derive the convergence rate of the incompressible limit, it suffices to deal with $H(t)$ and $K(t)$. First, recalling that $V=\divg u$ and $\dot V=\divg u_t+u\cdot\nabla\divg u$, by (\ref{2.1}), (\ref{2.4}), (\ref{2.5}), (\ref{3.2}) and (\ref{3.14}) we deduce that for any $2\leq r^*<\infty$,
\begin{equation}\label{4.39}
\begin{aligned}
\|u_t\|_{L^{r^*}}^2&\leq C\|\dot u\|_{L^{r^*}}^2+C\|u\|_{H^1}^2\|\nabla u\|_{H^1}^2\\
 & \leq C\left(\|\sqrt\rho \dot u\|_{L^2}^2+\|\nabla \dot u\|_{L^2}^2 +\|\nabla u\|_{H^1}^2\right)\\
 &\leq C\left(\|\sqrt\rho \dot u\|_{L^2}^2+\|\nabla^\bot\cdot \dot u\|_{L^2}^2 +\|\divg \dot u\|_{L^2}^2 +\|\nabla u\|_{H^1}^2\right)\\
  &\leq C\left(\|\sqrt\rho \dot u\|_{L^2}^2+\|\nabla^\bot\cdot \dot u\|_{L^2}^2 +\|\dot V\|_{L^2}^2+\|\nabla u\|_{L^4}^4 +\|\nabla u\|_{H^1}^2\right),
\end{aligned}
\end{equation}
and hence, it follows from (\ref{3.2}), (\ref{3.3}) and (\ref{4.1}) that
\begin{equation}\label{4.40}
\begin{aligned}
\int_0^T H(t)\md t&\leq C+C\int_0^T \left(\|\nabla u\|_{W^{1,q}}+\|\nabla u\|_{L^4}^2+ t\left\| \mathcal{P}u_t\right\|_{L^{r^*}}^2+ t\left\| \mathcal{Q}u_t\right\|_{L^{r^*}}^2\right)\md t\\
&\leq C+C\int_0^T  t \| u_t\|_{L^{r^*}}^2\md t\\
& \leq C+C\int_0^T t\left(\|\sqrt\rho \dot u\|_{L^2}^2+\|\nabla^\bot\cdot \dot u\|_{L^2}^2 +\|\dot V\|_{L^2}^2+\|\nabla u\|_{L^4}^4 +\|\nabla u\|_{H^1}^2\right)\md t\\
&\leq C.
\end{aligned}
\end{equation}

Thanks to (\ref{2.8}), (\ref{3.1}), (\ref{3.2}) and (\ref{4.1}), we have
\begin{equation}\label{4.41}
\begin{aligned}
&\nu^{-\frac{1}{2}}\int_0^T\left(\|\sqrt\rho u_t\|_{L^2}+\|\sqrt\eta v_t\|_{L^2}\right)\md t\\
&\quad\leq C\nu^{-\frac{1}{2}}\int_0^T\left(\|\sqrt\rho \dot u\|_{L^2}+\|u\|_{H^1}\|\nabla u\|_{H^1}+\|\sqrt\eta v_t\|_{L^2}\right)\md t\\
&\quad\leq C\nu^{-\frac{1}{2}}\int_0^T\left(1+\|\sqrt\rho \dot u\|_{L^2}+ \|\nabla u\|_{H^1}+\|\sqrt\eta v_t\|_{L^2}\right)\md t\\
&\quad\leq C\nu^{-\frac{1}{2}}\left(\int_0^T\left(1+\|\sqrt\rho \dot u\|_{L^2}^2+ \|\nabla u\|_{H^1}^2+\|\sqrt\eta v_t\|_{L^2}^2\right)\md t\right)^{\frac{1}{2}}\\
&\quad\leq C\nu^{-\frac{1}{2}},
\end{aligned}
\end{equation}
 
Analogously to the derivation of (\ref{4.40}), we easily infer from (\ref{4.39}) and H\"{o}lder inequality that
\begin{equation}\label{4.42}
\begin{aligned}
\nu^{-\frac{3}{4}}\int_0^T  \left(\sqrt t\|  u_t\|_{L^4}+\sqrt t\|  v_t\|_{L^4}\right)\md t &\leq C\nu^{-\frac{3}{4}}  \left(\int_0^T  t \left(\|  u_t\|_{L^4}^2+ \|  v_t\|_{L^4}^2\right)\md t\right)^{\frac{1}{2}}\\
&\leq C\nu^{-\frac{3}{4}},
  \end{aligned}
\end{equation}
and 
\begin{equation}\label{4.43}
\begin{aligned}
 \nu^{-\frac{2}{3}} \int_0^T\left(\sqrt t\|  u_t\|_{L^4}+\sqrt t\|  v_t\|_{L^4}\right)^{\frac{4}{3}}\md t &\leq C \nu^{-\frac{2}{3}}  \left(\int_0^T  t \left(\|  u_t\|_{L^4}^2+ \|  v_t\|_{L^4}^2\right)\md t\right)^{\frac{2}{3}}\\
 &\leq C\nu^{-\frac{2}{3}}.
 \end{aligned}
\end{equation}

Collecting (\ref{4.40})--(\ref{4.43}) together gives rise to
\begin{equation}\label{4.44}
\int_0^T K(t)\md t\leq C\left(\nu^{-1}+ \nu^{-\frac{1}{2}}+ \nu^{-\frac{3}{4}}+ \nu^{-\frac{2}{3}}\right)\leq C\nu^{-\frac{1}{2}}.
\end{equation}

It follows from   Cauchy-Schwarz inequality,  Poincar${\rm\acute{e}}$ inequality and (\ref{3.2}) that
\begin{align*}
&\left|2\int \eta (\mathcal{Q}u)  (\mathcal{P}u-v)\mdx +  \int \eta |\mathcal{Q}u|^2\mdx\right|\leq \frac{1}{4}\|\sqrt\eta(\mathcal{P}u-v)\|_{L^2}^2+C\|\mathcal{Q}u\|_{L^2}^2\\
&\quad \leq \frac{1}{4}\|\sqrt\eta(\mathcal{P}u-v)\|_{L^2}^2+C\|\divg(\mathcal{Q}u)\|_{L^2}^2\leq \frac{1}{4}\|\sqrt\eta(\mathcal{P}u-v)\|_{L^2}^2+C\nu^{-1},
\end{align*}
so that, by (\ref{4.40}) and  (\ref{4.44}) and Gronwall inequality, we conclude from (\ref{4.38}) that
\begin{equation}\label{4.45}
\begin{aligned}
& \sup_{0\leq t\leq T}\left(\|\sqrt\eta (\mathcal{P}u-v)(t)\|_{L^2}^2+\|\phi(t)\|_{L^2}^2+\frac{\|\phi(t)\|_{L^2}^2}{t}\right)\\
 &\qquad+\int_0^T\left( \|\nabla(\mathcal{P}u-v)\|_{L^2}^2+\frac{\|\phi\|_{L^2}^2}{t^2}\right)\md t\leq C \nu^{-\frac{1}{2}}.
 \end{aligned}
\end{equation}
This, combined with (\ref{4.11}) and (\ref{2.4}), shows
\begin{equation}\label{4.46}
\sup_{0\leq t\leq T}\|(\rho-\eta)(t)\|_{L^2}^2+\int_0^T \| \mathcal{P}u-v \|_{L^2}^2\md t\leq C \nu^{-\frac{1}{2}}.
\end{equation}

Finally, due to (\ref{1.10}), (\ref{3.2}) and (\ref{4.1}), it is easily seen that 
$$\|\nabla \mathcal{Q}u\|_{L^2}^2\leq C\|\divg u\|_{L^2}^2\leq C\nu^{-1}
$$ 
and
\begin{align*}
\int_0^T\|\nabla \mathcal{Q}u\|_{H^1}^2\md t&\leq C\int_0^T\| \divg u\|_{H^1}^2\md t\leq C\nu^{-2}\int_0^T\left(\|\rho\|_{H^1}^2+ \|F\|_{H^1}^2\right)\md t\\
&\leq C\nu^{-2}\int_0^T\left(1+ \|\nabla F\|_{L^2}^2\right)\md t\leq C\nu^{-2},
\end{align*}
which, together with (\ref{4.45}) and (\ref{4.46}), finishes the proof of Theorem \ref{thm1.2}.\hfill$\square$

\vskip 4mm

\noindent{\bf Data Availability.} Data sharing not applicable to this article as no datasets were
generated or analyzed during the current study.

\vskip 2mm

\noindent{\bf Declarations}
\vskip 2mm

\noindent{\bf Conflict of interest.} On behalf of all authors, the corresponding author states that there is no conflict of interest.

\end{document}